\documentclass[11pt,reqno]{amsart}

\usepackage{amsthm}
\usepackage{amssymb}
\usepackage{latexsym}
\usepackage{multicol}
\usepackage{verbatim,enumerate}
\usepackage{accents}
\usepackage[usenames]{color}
\usepackage[all, poly]{xy}
\usepackage[utf8]{inputenc}
\usepackage{hyperref}
\usepackage{amsmath, amscd}
\usepackage{soul}
\usepackage{tikz}
\usetikzlibrary{matrix,arrows,decorations.pathmorphing}
\usepackage{ytableau}

\theoremstyle{plain}
\newtheorem*{prop}{Proposition}
\newtheorem{thm}{Theorem}
\newtheorem*{lem}{Lemma}
\newtheorem*{cor}{Corollary}

\theoremstyle{definition}

\newtheorem*{defn}{Definition}

\newtheorem*{rem}{Remark}

\theoremstyle{remark}

\newcommand{\lie}[1]{\mathfrak{#1}}

\newcommand\bc{\mathbb C}

\newcommand{\Hom}{\operatorname{Hom}}

\def\gr{\operatorname{gr}}

\advance\textwidth by 1.2in  \advance\oddsidemargin by -.6in \advance\evensidemargin by -.6in

\newcounter{cnt}
 \makeatletter
\def\mydggeometry{\makeatletter\dg@YGRID=1\dg@XGRID=20\unitlength=0.003pt\makeatother}
\makeatother \theoremstyle{remark}

\numberwithin{equation}{section}
\makeatletter
\def\section{\def\@secnumfont{\mdseries}\@startsection{section}{1}%
  \z@{.7\linespacing\@plus\linespacing}{.5\linespacing}%
  {\normalfont\scshape\centering}}
\def\subsection{\def\@secnumfont{\bfseries}\@startsection{subsection}{2}%
  {\parindent}{.5\linespacing\@plus.7\linespacing}{-.5em}%
  {\normalfont\bfseries}}
\makeatother

\begin{document}

\title[Graded decompositions of fusion products in rank two]{Graded decompositions of fusion products in rank two}
\author{Leon Barth}
\address{Faculty of Mathematics Ruhr-University Bochum}
\email{leon.barth@rub.de}
\author{Deniz Kus}
\address{Faculty of Mathematics Ruhr-University Bochum}
\email{deniz.kus@rub.de}
\thanks{}

\begin{abstract}
We determine the graded decompositions of fusion products of finite-dimensional irreducible representations for simple Lie algebras of rank two. Moreover, we give generators and relations for these representations and obtain as a consequence that the Schur positivity conjecture holds in this case. The graded Littlewood-Richardson coefficients in the decomposition are parametrized by lattice points in convex polytopes and an explicit hyperplane description is given in the various types.
\end{abstract}

\maketitle


\section{Introduction}\label{section1}
Fusion products have been introduced in \cite{FL99} by proposing a way to construct generalized versions of Kostka polynomials. In turned out that these representations are very useful and generalize many families of representations, such as local Weyl modules, Kirillov-Reshetikhin modules and Demazure modules, which have been extensively studied due to their application to mathematical physics such as the X = M conjecture, and their rich combinatorial structure (see for example \cite{C01,CIK14,FK08,HKOTT02,Na11} and references therein). Recently, Naoi further highlighted the importance of fusion products by proving that the graded limit of a tensor product of (a certain class of) representations of quantum affine algebras is isomorphic to the fusion product of the corresponding representations for current algebras \cite{Na17}. Moreover, the characters of graded limits satisfy a graded version of the Q-system \cite{CV13,KV14} which is a certain difference equation arising from models in statistical mechanics \cite{HKOTT02}. Although the notion of fusion products is 20 years old, very little is known about their structure, especially in the case of finite-dimensional irreducible modules (only in special cases generators and relations are determined; see for example \cite{FF02,F15a,KL14}). The motivation of this paper is to settle this problem in rank two and to parametrize its highest weight vectors by lattice points in polytopes. \par

To each element $\lambda\in \mathbb{Z}_+^n$ we can associate an irreducible finite-dimensional representation of a complex finite-dimensional simple Lie algebra $\mathfrak{g}$ denoted by $V(\lambda)$. A very important issue originated from its applications to quantum physics is to understand their tensor product structure. The fusion product can be understood as a graded version of the tensor product and is denoted in the literature often by $V(\lambda)^{z_1}*V(\mu)^{z_2}$; its definition depends \textit{a priori} on a choice of distinct complex numbers $z_1,z_2\in\mathbb{C}$ (for details we refer to \cite{FL99}). The independence of the parameters was initially conjectured in \cite[Conjecture 1.8]{FL99} and is still an unsolved problem (except in the above mentioned cases where the structure has been determined). 

Another importance of fusion products is due to its connection to algebraic combinatorics. Schur functions $s_{\lambda}$ give a linear basis for the ring of symmetric functions and are labelled by partitions. Recently, a lot of work has gone into studying whether certain expressions of the form $s_{\lambda}s_{\mu}-s_{\rho}s_{\nu}$ are Schur positive, i.e. whether they can be written as a non-negative linear combination of Schur functions \cite{CFS12,FFLP05,LPP07}. The Schur positivity conjecture states under which assimptions on $\lambda,\mu,\rho,\nu$ the above expression is Schur positive and is known to be true only in some fairly simple special cases \cite{CFS12}. This conjecture can be generalized and stated in the language of Littlewood-Richardson coefficients (see \cite{CFS12}). As an application of the main theorem we will derive the Schur positivity conjecture in the rank two case (see Corollary~\ref{Schpo}).

We will now describe our results in more details. Denote by $\lie g[t]$ the algebra of polynomial maps $\mathbb{C}\rightarrow \mathfrak{g}$ and by $\mathcal{F}_{\lambda,\mu}$ the cyclic quotient of the universal algebra of $\lie g[t]$ generated by a non-zero vector $v$ subject to the defining relations described in Section~\ref{section31}.
Then, there is by construction a surjective homomorphism $\mathcal{F}_{\lambda,\mu}\rightarrow V(\lambda)*V(\mu)$ and it has been conjectured that this map is in fact an isomorphism. The conjecture is known to be true in the following special cases:
\begin{enumerate}
\item $\mathfrak{g}=\mathfrak{sl}_2$ \cite{FF02}
\item $\mathfrak{g}=\mathfrak{sl}_{n}$ and $\lambda,\mu$ are rectangular \cite{F15a,V13}
\item Arbitrary $\mathfrak{g}$ and $\lambda=\mu$ \cite{KL14}.
\end{enumerate}
We emphasize that $\mathfrak{sl}_2$ is the only case where the strucure has been fully understood. In this paper we prove the above conjecture for $\mathfrak{g}$ of rank two where we have to make the additional assumption in the case of $G_2$ that $\min\{\lambda(h_2),\mu(h_2)\}=0$ where $\alpha_1$ is the unique short root (see Theorem~\ref{mainthm}). We believe that this assumption is necessary to end up in a paramatrization of the graded Littlewood-Richardson coefficients by convex polytopes.

We outline the idea of the proof. Let $\lambda=m_1\varpi_1+m_2\varpi_2$ and $\mu=n_1\varpi_1+n_2\varpi_2$. We are interested in determining a parametrizing set of Littlewood-Richardson coefficients compatible with the fusion product grading. The strategy is to find an expression of the following form
$$V(\lambda)*V(\mu)=\sum_{\mathbf{s}\in \mathcal{S}}\mathbf{U}(\lie g\otimes 1)X_{\mathbf{s}}v,$$
where $\mathcal{S}\subseteq \mathbb{Z}_+^{|R_+|}$ is a finite subset parametrizing the highest weight vectors up to a filtration and $X_{\mathbf{s}}$ is defined in Section~\ref{section34}. To achieve such an expression, we will define suitable differential operators on enveloping algebras (similar as in \cite{BK15,FFL11, FFL2011}) and prove a straightening law (see Proposition~\ref{pr1}, Proposition~\ref{pr2}, and Proposition~\ref{pr3}).
It is striking that the proof of the straightening law for rank 2 Lie algebras leads to the fact that $\mathcal{S}$ can be described in terms of lattice points of a convex polytope (see Definition~\ref{maindefn}). Now comparing the number of lattice points with the Littlewood-Richardson coefficients (for example using the models in \cite{Li90} or \cite{PV05}) will give the desired result. \par

The paper is organized as follows. In Section~\ref{section2} we give the main notations. In Section~\ref{section3} we state the main result of the paper and give the hyperplane description of the corresponding polytopes in the various cases. Moreover, we discuss the strategy of the proof and its application to the Schur positivity conjecture. Section~\ref{section4} is devoted to the proof of the main theorem in type $A_2$,  Section~\ref{section5} deals with type $C_2$ and Section~\ref{section6} with type $G_2$.



\section{Notations}\label{section2}
In this section we recall the notion of fusion products \cite{FL99} and set up the notation needed in the rest of the paper.
\subsection{}Throughout this paper we denote by $\mathbb{C}$ the field of complex numbers and by $\mathbb{Z}$ (resp. $\mathbb{Z}_{+}$, $\mathbb{N}$) the subset of integers (resp. non-negative, positive integers). For $r\in\mathbb{Z}$ we set 
$$r_+=\max\{r,0\}.$$
Let $\lie g$  be a finite--dimensional complex simple Lie algebra and $\lie h\subseteq \mathfrak{g}$ a Cartan subalgebra  of $\lie g$. We denote the corresponding set of roots with respect to $\lie h$ by $R$ and fix $\Delta=\{\alpha_1,\dots,\alpha_n\}$ a basis for $R$. The corresponding sets of positive and negative roots are denoted as usual by $R^\pm$ and for $\alpha\in R$ we let $h_{\alpha}$ be the corresponding coroot of $\alpha$. We fix a Chevalley basis $\{x_{\alpha}^{\pm},\ h_i: \alpha\in R^+, 1\leq i\leq n\}$ of $\lie g$. Furthermore, let $Q$ be the root lattice and $Q^{+}$ the subset of non-negative integer linear combinations of simple roots. Let $P$ be the weight lattice of $\mathfrak{g}$ with basis $\{\varpi_1,\dots,\varpi_n\}$ and $P^+$ be the subset of dominant integral weights.
\subsection{}We let $\mathbf U(\lie g)$ be the universal enveloping algebra of $\lie g$ and $\lie g[t]$ be the Lie algebra of polynomial maps from $\mathbb C$ to $\lie g$ with the obvious pointwise Lie bracket: $$[x\otimes f, y\otimes g]=[x,y]\otimes fg,\ \ x,y\in\lie g,\ \ f,g\in\mathbb C[t].$$ The Lie algebra $\lie g[t]$  and its universal enveloping algebra inherit a grading from the degree grading of $\mathbb C[t]$. Thus an element $(a_1\otimes t^{r_1})\cdots (a_s\otimes t^{r_s})$, $a_j\in\lie g$, $r_j\in\mathbb Z_+$ for $1\le j\le s$ will have grade $r_1+\cdots+r_s$. We denote by $\mathbf U(\lie g[t])_k$ be the homogeneous component of degree $k$ and recall that it is a $\lie g$--module for all $k\in\mathbb{Z}_+$. Suppose now that we are given a cyclic $\lie g[t]$--module $V$ generated by a vector $v$. Define an increasing filtration $0\subseteq V^0\subseteq V^1\subseteq\cdots $ of $\lie g$--submodules of $V$ by $$V^k=\sum_{s=0}^k \mathbf U(\lie g[t])_s v.$$ The associated graded vector space $\gr V$ admits an action of $\lie g[t]$ given by: 
$$
x(v+V^{k})= x.v+ V^{k+s},\ \ x\in\lie g[t]_s,\ \ v\in V^{k+1}.
$$
Furthermore, $\gr V$ is a cyclic $\lie g[t]$--module with cyclic generator $\bar v$, the image of $v$ in $\gr  V$. The fusion product \cite{FL99} is a $\lie g[t]$--module of the form $\gr V$ for a special choice of cyclic $\lie g[t]$--module $V$ which we define next. Given a dominant integral weight $\lambda\in P^+$ and $z\in\mathbb \bc$, let $V(\lambda)^z$ be the $\lie g[t]$--module with action 
$$(x\otimes t^r)w=z^rx.w,\ x\in \lie g,\ w\in V(\lambda),\ r\in \mathbb Z_+.$$
Given $\lambda,\mu\in P^+$ and distinct complex numbers $z_1,z_2$, the Chinese remainder theorem implies that the tensor product $V_{\lambda,\mu}^{z_1,z_2}:=V(\lambda)^{z_1}\otimes V(\mu)^{z_2}$ is cyclic and generated by $v_1\otimes v_2$ where $v_1$ and $v_2$ respectively is the highest weight vector of $V(\lambda)$ and $V(\mu)$ respectively.
The fusion product is defined as follows
$$V(\lambda)^{z_1}*V(\mu)^{z_2}:=\gr (V_{\lambda,\mu}^{z_1,z_2}).$$ 
\begin{rem}The fusion product can be defined more generally for any collection of dominant integral weights $\lambda_1,\dots,\lambda_N\in P^+$ and a tuple $(z_1,\dots,z_N)$ of pairwise distinct complex numbers (see \cite{FL99} for details). Clearly the definiton depends on the parameters $z_i$, $1\le i\le N$. However, it is conjectured in \cite{FL99} (and proved in several special cases) that the structure is independent of the choice. If $N=2$, the independence can be proven directly. 
\end{rem}
\begin{lem}\label{inde} For each pair $(z_1,z_2)$ of distinct complex numbers, we have an isomorphism of $\mathfrak{g}[t]$-modules
$$V(\lambda)^{z_1}*V(\mu)^{z_2}\rightarrow V(\lambda)^{0}*V(\mu)^{1}.$$
\begin{proof}
In order to distinguish the highest weight elements, we shall write $v_1^{z_1}\otimes v_2^{z_2}$ and $v_1^{0}\otimes v_2^{1}$ respectively. To see the isomorphism, we only have to observe that
$$X(\overline{v_1^{z_1}\otimes v_2^{z_2}})=0\Rightarrow X(\overline{v_1^{0}\otimes v_2^{1}})=0,\ \ X\in \mathbf{U}(\mathfrak{g}[t]).$$
\textit{Claim:} For each homogeneous element $X\in \mathbf{U}(\mathfrak{g}[t])$   of degree $k$ there exists suitable elements $Y_i$ of degree strictly smaller than $k$ such that \begin{equation}\label{cl1}X(v_1^{z_1}\otimes v_2^{z_2})=(z_2-z_1)^k X(v_1^{0}\otimes v_2^{1})+\sum_{i}Y_i (v_1^{z_1}\otimes v_2^{z_2}).\end{equation}
\textit{Proof of the Claim:} It is enough to show the claim for monomials $X$. Since $(x\otimes t^r)(\overline{v_1^{z_1}\otimes v_2^{z_2}})=0$ for all $r\geq 2$, we can additionally assume by the PBW theorem that $X$ is of the form
$$X=(x_1\otimes 1)\cdots (x_r\otimes 1)(y_1\otimes t)\cdots (y_k\otimes t).$$
Modulo terms of strictly smaller degree (applied to $v_1^{z_1}\otimes v_2^{z_2}$) we obtain
\begin{align*}X(v_1^{z_1}\otimes v_2^{z_2})&=(x_1\otimes 1)\cdots (x_r\otimes 1)(y_1\otimes t)\cdots (y_{k-1}\otimes t)(z_1 y_kv_1^{z_1}\otimes v_2^{z_2}+z_2 v_1^{z_1}\otimes y_kv_2^{z_2})&\\&=(z_2-z_1)(x_1\otimes 1)\cdots (x_r\otimes 1)(y_1\otimes t)\cdots (y_{k-1}\otimes t)(v_1^{z_1}\otimes y_kv_2^{z_2})&\\&=\cdots =(z_2-z_1)^{k}(x_1\otimes 1)\cdots (x_r\otimes 1)(v_1^{z_1}\otimes y_1\cdots y_kv_2^{z_2})&\\&=(z_2-z_1)^{k}X(v^0_1\otimes v_2^1).
\end{align*}
This shows the claim. \par
Now if $X(\overline{v_1^{z_1}\otimes v_2^{z_2}})=0$, i.e. $X(v_1^{z_1}\otimes v_2^{z_2})$ can be written as a sum of elements of strictly smaller degree (applied to $v_1^{z_1}\otimes v_2^{z_2}$), we obtain from \eqref{cl1} that
$$X(v_1^{0}\otimes v_2^{1})=\sum_{i}Z_i (v_1^{z_1}\otimes v_2^{z_2})$$
for suitable elements $Z_i$ of degree strictly smaller than $k$. In order to show that 
$X(v_1^{0}\otimes v_2^{1})$ can be written as a sum of elements of strictly smaller degree than $k$ (applied to $v_1^{0}\otimes v_2^{1}$) we use once more \eqref{cl1} and replace each element $Z_i (v_1^{z_1}\otimes v_2^{z_2})$ in the above sum. This gives the desired result since $Y(v_1^{z_1}\otimes v_2^{z_2})=(v_1^{0}\otimes v_2^{1})$ for $Y\in\mathbf{U}(\lie g)$.
\end{proof}
\end{lem}
In the rest of the paper we simply write $V(\lambda)*V(\mu)$ since the fusion product is independent of the choice of parameters by Lemma~\ref{inde}.
Although fusion products have been introduced around 20 years ago, very little is known about their structure; only in special cases generators and relations are determined (see for example \cite{F15a,KL14}) and the complete structure has only been understood in the case of $A_1$ \cite{FF02}. The motivation of this paper is to investigate the structure of fusion products in rank two.
\section{The Main Theorem}\label{section3}
In this section we state our main theorem and discuss the strategy of the proof as well as several consequences.
\subsection{}\label{section31} We denote by $\mathcal{F}_{\lambda,\mu}$ the cyclic $\mathbf{U}(\lie g[t])$-module generated by a non-zero vector $v$ subject to the relations 
$$\lie n^+[t]v=0,\ (h\otimes t^r)v=\delta_{r,0}(\lambda+\mu)(h)v,\ \forall h\in\lie h$$
$$(x^{-}_{\alpha}\otimes 1)^{(\lambda+\mu)(h_{\alpha})+1}v=(x^{-}_{\alpha}\otimes t)^{{\text{min}\{\lambda(h_{\alpha}),\mu(h_{\alpha})\}}+1}v=0,\ \forall \alpha\in R^+$$
$$(x^{-}_{\alpha}\otimes t^r)v=0,\ \forall \alpha\in R^+,\ r\geq 2.$$
Obviously, there is a surjective homomorphism 
\begin{align}\label{map31}
	\mathcal{F}_{\lambda,\mu}\rightarrow V(\lambda)*V(\mu)
\end{align}
 and it is still an open problem, whether the above map is an isomorphism. 
\subsection{} In this subsection we will state our main theorem.
\begin{defn}\label{maindefn}Let $\lie g$ be of rank two and $\lambda,\mu\in P$. We set $m_i=\lambda(h_{i}),\ n_i=\mu(h_{i})$ for $i=1,2$.
\begin{enumerate}
\item Let $\mathcal{P}_{\lambda,\mu}^{A}\subseteq \mathbb{R}_+^3$ be the convex polytope defined by the inequalities 
\begin{gather}
	\label{1}a\leq \min\{m_1,n_1\},\ c\leq \min\{m_2,n_2\},\ a+b+c\leq \min\{m_1+m_2,n_1+n_2\},\\
	\label{2}2a+b\leq m_1+n_1,\ \ 2c+b\leq m_2+n_2.
\end{gather}

\item Let $\mathcal{P}_{\lambda,\mu}^{C}\subseteq \mathbb{R}_+^4$ be the convex polytope defined by the inequalities 
\begin{gather}
	\label{3}a\leq \min\{m_1,n_1\},\ d\leq \min\{m_2,n_2\},\\ 
	\label{3.1}  a+b+c\leq \min\{m_1+m_2,n_1+n_2\},\ a+b+d\leq \min\{m_1+m_2,n_1+n_2\},\\
	\label{4}2a+b\leq m_1+n_1,\ \ 2d+b\leq m_2+n_2,\\
	\label{5}2a+b+2(c-d)\leq m_1+n_1.
\end{gather}
\item Let $\mathcal{P}_{\lambda,\mu}^{G}\subseteq \mathbb{R}_+^6$ be the convex polytope defined by the inequalities 
\begin{gather}
	\label{6}a\leq \min\{m_1,n_1\},\ \ b\leq m_2+n_2,\ \ f\leq \min\{m_2,n_2\},\ \ b+e-a\leq m_2+n_2,\\ 
	\label{7} a+c+d\leq \min\{m_1+m_2,n_1+n_2\},\ \ a+b+c\leq \min\{m_1+m_2,n_1+n_2\},\\
	\begin{split}
		\label{71} a+b+c+d\leq \min\{m_1+2m_2,n_1+2n_2\},\\ b+c+d+e\leq \min\{m_1+2m_2,n_1+2n_2\},	
	\end{split}\\
	\label{8}2(a+c)+3d-b\leq m_1+n_1,\ \ 2(a+c)+b+d\leq m_1+n_1.
\end{gather}
\end{enumerate}
\end{defn}
The upper index denotes the type of the Lie algebra. Note that $\mathcal{P}_{\lambda,\mu}^{\lie g}$ is empty if $\lambda$ or $\mu$ is not a dominant integral weight. We denote the lattice points in $\mathcal{P}_{\lambda,\mu}^{\lie g}$ by $\mathcal{S}_{\lambda,\mu}^{\lie g}$ and for $\mathbf{s}\in \mathcal{S}_{\lambda,\mu}^{\lie g}$ we define its degree $\deg(\mathbf{s})$ as the sum over its entries and its weight by $$\mathrm{wt}(\mathbf{s})=\begin{cases}(a+b)\alpha_1+(b+c)\alpha_2,& \text{if $\lie g=A_2$},\\ (a+b+2c)\alpha_1+(b+c+d)\alpha_2,& \text{if $\lie g=C_2$},\\ (a+b+2c+3d+3e)\alpha_1+(b+c+d+2e)\alpha_2,& \text{if $\lie g=G_2$}. \end{cases}$$
In our convention $\alpha_1$ is the unique short root in the case of $C_2$ and $G_2$.
\begin{thm}\label{mainthm} Let $\mathfrak{g}$ be of rank two, $\lambda,\mu\in P^+$ be dominant integral weights such that $\lambda$ or $\mu$ is a multiple of a minuscule fundamental weight in the case of $G_2$. 
\begin{enumerate}
\item We have an isomorphism of $\lie g[t]$-modules 
$$\mathcal{F}_{\lambda,\mu}\xrightarrow{\sim} V(\lambda)*V(\mu).$$
\item The multiplicity of $V(\nu)$ in degree $r\in \mathbb{Z}_+$ in the fusion product $V(\lambda)*V(\mu)$ is given by the cardinality of 
$$\{\mathbf{s}\in \mathcal{S}_{\lambda,\mu}^{\lie g}: \mathrm{wt}(\mathbf{s})=\lambda+\mu-\nu,\ \ \deg(\mathbf{s})=r\}.$$
\end{enumerate}
\end{thm}
\begin{rem}\label{remrec}We see that each lattice point in the polytope $\mathcal{P}_{\lambda,\mu}^{\lie g}$ is also a lattice point in the corresponding FFLV/Gornitskii polytope $P(\lambda)$ and $P(\mu)$ respectively (see \cite{FFL11, FFL2011,G15}). The restriction on the dominant integral weights in the theorem for type $G_2$ is necessary and preliminary calculations show that our methods do not lead to a convex polytope. We expect a union of polytopes in the general case for type $G_2$.
\end{rem}
\subsection{} We first discuss a consequence of our main theorem. Schur functions $s_{\lambda}$ give a linear basis for the ring of symmetric functions and are labelled by partitions. Recently, a lot of work has gone into studying whether certain expressions of the form $s_{\lambda}s_{\mu}-s_{\rho}s_{\nu}$ are Schur positive, i.e. whether they can be written as a non-negative linear combination of Schur functions. Since Schur functions describe the characters of irreducible finite-dimensional modules in type $A$, the above question has been generalized and formulated in terms of Littlewood Richardson coefficients in \cite[Conjecture 2.3]{CFS12}. So far the aforementioned conjecture is known to be true only in some fairly simple special cases \cite{CFS12}, namely for multiples of fundamental minuscule weights or if $\mathfrak{g}$ is of type $A_2$. 
\begin{cor}\label{Schpo} Let $\mathfrak{g}$ be of rank two and $\lambda_i,\mu_i\in P^+$, $i=1,2$ be dominant integral weights such that $\lambda_1+\lambda_2=\mu_1+\mu_2$. In the case of $G_2$ we additionally assume that $\lambda_2$ and $\mu_2$ are multiples of a fundamental miniscule weight.  If $\min\{\lambda_1(h_{\alpha}),\lambda_2(h_{\alpha})\}\leq \min\{\mu_1(h_{\alpha}),\mu_2(h_{\alpha})\} \ \forall \alpha\in R^+$, we have
\begin{equation*}
\dim\Hom_{\lie g}  (V(\nu), V(\lambda_1) \otimes V(\lambda_2)) \le \dim\Hom_{\lie g}(V(\nu), V(\mu_1) \otimes V(\mu_2)  ), \text{ for all $\nu \in P^+$}.
\end{equation*}
\begin{proof}
Since $V(\lambda_1)*V(\lambda_2)\cong V(\lambda_1)\otimes V(\lambda_2)$ as a $\mathfrak{g}$-module, we only have to prove that there exists a surjective map $V(\mu_1)*V(\mu_2)\rightarrow V(\lambda_1)*V(\lambda_2)$. Since $\mathcal{F}_{\lambda_1,\lambda_2}\cong V(\lambda_1)*V(\lambda_2)$, this is a direct consequence of Theorem~\ref{mainthm}(1) and the assumptions.
\end{proof}
\end{cor}
With other words, we can extend \cite[Conjecture 2.3]{CFS12} to the simple Lie algebra of type $C_2$ and $G_2$ (under the above restriction) and provide a completely different proof in the case of $A_2$.
\subsection{}\label{section34} In this subsection we explain the strategy of the proof of Theorem~\ref{mainthm}. Note that the first part of Theorem~\ref{mainthm} would follow from the injectivity of  the map in \eqref{map31}. This again is a consequence of the multiplicity bound
\begin{equation}\label{multb}\sum_{\nu\in P^+}[\mathcal{F}_{\lambda,\mu}: V(\nu)]\leq \sum_{\nu\in P^+} [V(\lambda)\otimes V(\mu): V(\nu)].\end{equation}
The strategy in showing such an inequality \eqref{multb} will be as follows. We write in the rest of the paper for a positive root $\alpha=r_1\alpha_1+r_2\alpha_2$ simply $\alpha=(r_1,r_2)$. We define for an arbitrary tuple $\mathbf{s}$ of non-negative integers the monomial
\begin{align*}
	X_{\mathbf{s}}\!=\!\big(x^-_{(1,0)}\!\otimes t\big)^{\!(a)}\!\big(x^-_{(1,1)}\!\otimes t\big)^{\!(b)} 
	\!\times\!
	\begin{cases}
		\!\big(x^-_{(0,1)}\!\otimes t\big)^{\!(c)}&\!\! \!\text{if $\lie g=A_2$},\\ 
		\!\big(x^-_{(2,1)}\!\otimes t\big)^{\!(c)}\!\big(x^-_{(0,1)}\!\otimes t\big)^{\!(d)}& \!\!\!\text{if $\lie g=C_2$},\\ 
		\!\big(x^-_{(2,1)}\!\otimes t\big)^{\!(c)}\!\big(x^-_{(3,1)}\!\otimes t\big)^{\!(d)}\!\big(x^-_{(3,2)}\!\otimes t\big)^{\!(e)}\!\big(x^-_{(0,1)}\!\otimes t\big)^{\!(f)}& \!\!\!\text{if $\lie g=G_2$}. 
\end{cases}	
\end{align*}
where $x^{(s)}$ denotes the divided powers, i.e. $x^{(s)}=\frac{x^s}{s!}$.

\begin{prop}\label{wichtp}Under the assumptions of Theorem~\ref{mainthm} we have $$\mathcal{F}_{\lambda,\mu}=\sum_{\mathbf{s}\in \mathcal{S}_{\lambda,\mu}^{\lie g}}\mathbf{U}(\lie g\otimes 1)X_{\mathbf{s}}v,\ \ \ |\mathcal{S}_{\lambda,\mu}^{\lie g}|=\sum_{\nu\in P^+} [V(\lambda)\otimes V(\mu): V(\nu)].$$
Moreover, we can order the elements in $\mathcal{S}_{\lambda,\mu}^{\lie g}$ in a way such that $$(\mathfrak{n}^+\otimes 1)X_{\mathbf{s}}v\in \mathrm{span}\left\{X_{\mathbf{s}'}v: \mathbf{s}' \prec \mathbf{s}'\right\}.$$
\end{prop}
The above proposition implies \eqref{multb} and hence the first part of our main theorem. To see this, consider the finite set $\mathcal{S}_{\lambda,\mu}^{\lie g}$ and the aforementioned order $\mathbf{s}_1\prec \cdots \prec \mathbf{s}_r$ (which we will define later) on the elements of $\mathcal{S}_{\lambda,\mu}^{\lie g}$.
Define $$F_i=\sum_{j\leq i}\mathbf{U}(\lie g\otimes 1)X_{\mathbf{s}_j}v$$
which gives an increasing filtration 
$$0\subseteq F_1\subseteq F_2\subseteq \cdots \subseteq F_{r}=\mathcal{F}_{\lambda,\mu}.$$
If we act with $(\mathfrak{n}^+\otimes 1)$ on a monomial $X_{\mathbf{s}_i}v$, we obtain an element which is in the span of monomials corresponding to strictly smaller elements, i.e. $\mathbf{s}_1,\dots,\mathbf{s}_{i-1}$. This implies that we have a well--defined surjective map of $\lie g$--modules $V(\lambda+\mu-\text{wt}(\mathbf{s}_i))\rightarrow F_i/F_{i-1}$, which gives
$$|\mathcal{S}_{\lambda,\mu}^{\lie g}|\geq \sum_{\nu\in P^+}[\mathcal{F}_{\lambda,\mu}: V(\nu)].$$
Together with Proposition~\ref{wichtp} we get \eqref{multb}. The second part of the theorem follows by construction. Hence we are left to show Proposition~\ref{wichtp} which is done in Section~\ref{section4} for type $A_2$, Section~\ref{section5} for type $C_2$ and Section \ref{section6} for type $G_2$. 

\section{The case $A_2$}\label{section4}
In this section we prove Proposition~\ref{wichtp} for the rank two Lie algebra of type $A_2$.
\subsection{}We consider the following total order on $\mathbb{Z}^3_+$. We say that $\mathbf{s}=(a,b,c)\succ \mathbf{s}'=(a',b',c')$ if and only if 
one of the following conditions hold
\begin{center}
$\circ \  c<c'\ \ \ \quad \circ \ c=c',\  a<a'\ \ \ \quad 
\circ\  c=c',\ a=a',\ b<b'.$ 
\end{center}
\begin{prop}\label{pr1}Assume that $\mathbf{s}\in \mathbb{Z}_+^3$ is such that $\mathbf{s}\notin \mathcal{S}_{\lambda,\mu}^{A}$. Then we have
$$\left(X_{\mathbf{s}}+\sum_{\mathbf{s}'\prec\mathbf{s} } c_{\mathbf{s}'} X_{\mathbf{s}'}\right)v=0.$$
for suitable elements $c_{\mathbf{s}'}\in\mathbf{U}(\mathfrak{n}^-)$.
\begin{proof}
If $\mathbf{s}$ violates one of the inequalities in \eqref{1}, we would have 
\begin{equation}\label{wasw}\big(x^{-}_{(1,0)}\otimes t\big)^{(a)}\big(x^{-}_{(1,1)}\otimes t\big)^{(b)}\big(x^{-}_{(0,1)}\otimes t\big)^{(c)}v=0\end{equation}
and the claim is obvious. To see this, assume that $a+b+c>\min\{m_1+m_2,n_1+n_2\}$, i.e. 
\begin{equation}\label{wasw1}\big(x^{+}_{(1,0)}\otimes 1\big)^{(c)}\big(x^{+}_{(0,1)}\otimes 1\big)^{(a)}\big(x^{-}_{(1,1)}\otimes t\big)^{(a+b+c)}v=0\end{equation} by the relations of the module $\mathcal{F}_{\lambda,\mu}$.
Since $\mathfrak{n}^+[t]v=0 =(\mathfrak{h}\otimes t^r)v$ for all $r\geq 1$,
we obtain that equation \eqref{wasw1} is equal to \eqref{wasw}. If one of the other two inequalites in \eqref{1} are violated, we would immeaditely get \eqref{wasw}. So assume that $\mathbf{s}$ violates the first inequality in \eqref{2}. We consider 
$$0=\big(x^{-}_{(0,1)}\otimes t\big)^{(b+c)}\big(x^{+}_{(1,0)}\otimes t\big)^{(a)}\big(x^{-}_{(1,0)}\otimes 1\big)^{(2a+b)}v.$$
Expanding the above product we obtain up to sign 
\begin{align*}
	0&=\big(x^{-}_{(0,1)}\otimes t\big)^{(b+c)}\big(x^{-}_{(1,0)}\otimes 1\big)^{(b)}\big(x^{-}_{(1,0)}\otimes t\big)^{(a)}v&\\
	&=\sum_{i=0}^{b}\pm \big(x^{-}_{(1,0)}\otimes 1\big)^{(b-i)}\big(x^{-}_{(0,1)}\otimes t\big)^{(b+c-i)}\big(x^{-}_{(1,1)}\otimes t\big)^{(i)}\big(x^{-}_{(1,0)}\otimes t\big)^{(a)}v.
\end{align*}
The maximal term in the above sum is given if we minimalize the power of $\big(x^{-}_{(0,1)}\otimes t\big)$, so if $i=b$. This shows the claim. If $\mathbf{s}$ violates the second inequality in \eqref{2}, then consider
$$0=\big(x^{-}_{(1,0)}\otimes t\big)^{(a+b)}\big(x^{+}_{(0,1)}\otimes t\big)^{(c)}\big(x^{-}_{(0,1)}\otimes t\big)^{(2c+b)}v.$$
The proof is exactly the same and will be omitted. 
\end{proof}
\end{prop}
Note that Proposition~\ref{pr1} finishes all parts of Proposition~\ref{wichtp} except the statement on the cardinality of $\mathcal{S}_{\lambda,\mu}^{A}$.
\subsection{} There are several well-known formulas for the Littlewood-Richardson coefficients; however to the best of our knowledge none of them is compatible with the fusion product grading.  Consider the set $\mathcal{T}^{A}_{\lambda,\mu}$ consisting of triples $(a,b,c) \in \mathbb{Z}_{+}^3$ satisfying the following inequalities: 
	\begin{gather*}
		b  \le \min\{m_2, n_1\}, \qquad a+b-c \le n_1, \qquad b+c-a  \le m_2, \label{eq:pak1} \\
		c  \le n_2 , \qquad a  \le m_1,  \label{eq:pak2}\\
		2a+b-c  \le m_1 + n_1 , \qquad 2c+b-a  \le m_2 + n_2. \label{eq:pak3}
	\end{gather*}	
Then it has been proved in \cite[Section 3]{PV05} that 
$$|\mathcal{T}^{A}_{\lambda,\mu}|=\sum_{\nu\in P^+}[V(\lambda)\otimes V(\mu): V(\nu)].$$
So the following lemma finishes the proof of Proposition~\ref{wichtp} in the case of $A_2$.
\begin{lem}\label{}
There exists a bijective map $\mathcal{S}^{A}_{\lambda,\mu} \rightarrow  \mathcal{T}^{A}_{\lambda,\mu}$.
\begin{proof} 
First we define an equivalence relation on $\mathcal{S}^{A}_{\lambda,\mu} \cup \mathcal{T}^{A}_{\lambda,\mu}$ as follows
	\begin{align*}
		(a,b,c) \sim (a',b',c') : \Leftrightarrow \exists \ell \in \mathbb{Z} \colon (a', b' ,c')=(a+\ell,b-\ell,c+\ell).
	\end{align*}
Then $\mathcal{S}^{A}_{\lambda,\mu} \cup \mathcal{T}^{A}_{\lambda,\mu}$ decomposes into the disjoint union of equivalence classes. The following statements are straightforward to check:
		\begin{enumerate}
			\item[(i)] $(a,b,c) \in \mathcal{S}^{A}_{\lambda,\mu}  \text{ and }a,c \ge 1 \Rightarrow (a-1,b+1,c-1) \in \mathcal{S}^{A}_{\lambda,\mu}  $ \vspace{0,2cm}
		\item [(ii)]	$(a,b,c), (a+\ell,b-\ell,c+\ell)  \in \mathcal{T}^{A}_{\lambda,\mu} \text{ with } 1 \le \ell \le b  \Rightarrow (a+1, b-1, c+1)  \in \mathcal{T}^{A}_{\lambda,\mu}$
		\end{enumerate}
Now the existence of a bijective map $\mathcal{S}^{A}_{\lambda,\mu}  \rightarrow \mathcal{T}^{A}_{\lambda,\mu} $ follows if we can show that
for any equivalence class $M$ the cardinality of $M \cap \mathcal{S}^{A}_{\lambda,\mu}$ is equal to the cardinality of $M \cap \mathcal{T}^{A}_{\lambda,\mu} $. If $M \cap \mathcal{S}^{A}_{\lambda,\mu} $ is empty, there is nothing to show. So assume that the intersection is non-empty. By definition of the equivalence relation and (i) we can always choose an element $(a,b,c)\in M \cap \mathcal{S}^{A}_{\lambda,\mu}$ satisfying $a=0$ or $c=0$. Since the proof in both cases is similar, we only consider the case where $M$ is an equivalence class and $M \cap \mathcal{S}^{A}_{\lambda,\mu}$ contains an element of the form $(0,b,c)$. We fix this element in the rest of the proof and set 
\begin{align*}
			r & = \min \{ b, m_1, n_2 - c \}, \\ 
			R_1 & = \max \{ b-\min \{m_2,n_1\}, b-c -n_1, b+c  - m_2 \}, \\ 
		R_2 & = \min\big\{\min\{m_1,n_1\}, \min\{m_2,n_2\}-c,\min\{m_1+m_2,n_1+n_2\}-b-c,\\ & \quad\quad\quad \quad m_1+n_1-b, m_2+n_2-2c-b,b\big\}.
		\end{align*}
From the definition of the polytopes one quickly observes that (using (i) and (ii))
		$$
			(\ell,b-\ell,c+\ell) \in \mathcal{S}^{A}_{\lambda,\mu}  \Leftrightarrow 0\leq \ell \le R_2 \ \text{ and } \  (\ell,b-\ell,c+\ell) \in \mathcal{T}^{A}_{\lambda,\mu}  \Leftrightarrow R_1 \le \ell \le r.
		$$
So in order to finish the proof we are left to show that $R_2+1 = r-R_1+1$.
Note that \begin{align*} 
				R_1 & = \max \lbrace b-\min \{m_2,n_1\}, b-c -n_1, b+c  - m_2 \rbrace
				= \max \lbrace b-\min \{m_2,n_1\}, b+c  - m_2 \rbrace \\ 
				&= - \min \lbrace  \min \{m_2,n_1\} -b , m_2 - b - c \rbrace
				 = -\min \lbrace n_1- b, m_2 - b - c \rbrace \\
				&= b -\min \lbrace n_1, m_2 - c \rbrace.
			\end{align*}
Hence we get
\begin{align*}\min\{ r, r-R_1\}&=\min\big\{\min \lbrace b, m_1, n_2 - c \rbrace, \min \{ b, m_1, n_2 - c \rbrace - b + \min\{n_1, m_2 - c\}\big\}&\\& = \min\big\{\min \lbrace b, m_1, n_2 - c \rbrace, n_1, m_1+n_1-b, n_1+n_2 - b-c, &\\& \quad \quad \quad \quad \quad \quad m_2-c,m_1+m_2-b-c,m_2+n_2-b-2c \big\}=R_2
\end{align*}
which finishes the proof.
\end{proof}
\end{lem}
\section{The case $C_2$}\label{section5}
In this section we prove Proposition~\ref{wichtp} for the rank two Lie algebra of type $C_2$.
\subsection{}Again we define a total order on $\mathbb{Z}^4_+$ as follows. 
We say $\mathbf{s}=(a,b,c,d)\succ \mathbf{s}'=(a',b',c',d')$ if and only if one of the following conditions hold
\begin{center}
$\circ\ d<d',\ \  \quad \circ\ d=d',\ a<a',\  \ \quad \circ\ d=d',\ a=a',\ b<b',\ \  \quad \circ\ d=d',\ a=a',\ b=b',\ c<c'$.
\end{center}
\begin{prop}\label{pr2}
Assume that $\mathbf{s}\in \mathbb{Z}_+^4$ is such that $\mathbf{s}\notin \mathcal{S}^C_{\lambda,\mu}$ Then we have
$$\left(X_{\mathbf{s}}+\sum_{\mathbf{s}'\prec\mathbf{s} } c_{\mathbf{s}'} X_{\mathbf{s}'}\right)v=0.$$
for suitable elements $c_{\mathbf{s}'}\in\mathbf{U}(\mathfrak{n}^-)$.
\begin{proof}
If $(a,b,c,d)$ violates one of the inequalities in \eqref{3}, then the left hand side is already zero by the defining equalities of $\mathcal{F}_{\lambda,\mu}$. If the element violates the first inequality in \eqref{3.1}, then
\begin{align*}0&=\big(x^{+}_{(1,0)}\otimes t\big)^{(b)}\big(x^{+}_{(1,1)}\otimes t\big)^{(a)}\big(x^{-}_{(2,1)}\otimes 1\big)^{(a+b+c)}v&\\&=\big(x^{+}_{(1,0)}\otimes t\big)^{(b)}\big(x^{-}_{(1,0)}\otimes t\big)^{(a)}\big(x^{-}_{(2,1)}\otimes 1\big)^{(b+c)}v&\\&=\big(x^{-}_{(1,0)}\otimes t\big)^{(a)}\big(x^{+}_{(1,0)}\otimes t\big)^{(b)}\big(x^{-}_{(2,1)}\otimes 1\big)^{(b+c)}v
&\\&=(x^{-}_{(1,0)}\otimes t)^{(a)}(x^{-}_{(1,1)}\otimes t)^{(b)}(x^{-}_{(2,1)}\otimes 1)^{(c)}v\end{align*}
Now multiplying the above equation by $\big(x_{(0,1)}^{-}\otimes t\big)^{(d)}$ gives the claim.
The proof for the case when $(a,b,c,d)$ violates the second inequality in \eqref{3.1} is similar and will be omitted.
So assume now that the element violates the first inequality in \eqref{4}; the second is again similarly proven and will be omitted. Then we consider (again up to sign)
$$0=\big(x_{(0,1)}^-\otimes t\big)^{(b+d)}\big(x_{(1,0)}^+\otimes t\big)^{a}\big(x_{(1,0)}^-\otimes 1\big)^{(2a+b)}v=\big(x_{(1,0)}^-\otimes t\big)^{(a)}\big(x_{(0,1)}^-\otimes t\big)^{(b+d)}\big(x_{(1,0)}^-\otimes 1\big)^{(b)}v.$$
Now passing successively $\big(x_{(0,1)}^{-}\otimes t\big)$ through $\big(x_{(1,0)}^-\otimes 1\big)$ gives (using \cite[Lemma 4.1]{SaCh})
$$0=\sum \pm \big(x^{-}_{(1,0)}\otimes 1\big)^{(b-i-2j)}\big(x^{-}_{(1,0)}\otimes t\big)^{(a)}\big(x^{-}_{(1,1)}\otimes t\big)^{(i)}(x^{-}_{(2,1)}\otimes t\big)^{(j)}\big(x^{-}_{(0,1)}\otimes t\big)^{(b+d-i-j)}v,$$
where the sum runs over tuples $(i,j)\in \mathbb{Z}_+^2$ satisfying $i+j\leq b+d$ and $i+2j\leq b$. Hence the maximal monomial is given by $i=b$ and $j=0$. Now multiplying the above equation with $\big(x^{-}_{(2,1)}\otimes t\big)^{(c)}$ gives the result in this case.
If $(a,b,c,d)$ violates the inequality in \eqref{5}, then
$$\big(x^{-}_{(1,0)}\otimes 1\big)^{(2a+b+2c)}\big(x^{-}_{(0,1)}\otimes t\big)^{(d)}v=0$$ and hence it follows up to sign that 
\begin{align*}
0&=\big(x^{-}_{(0,1)}\otimes t\big)^{(b+c)}\big(x^{+}_{(1,0)}\otimes t\big)^{(a)}\big(x^{-}_{(1,0)}\otimes 1\big)^{(2a+b+2c)}\big(x^{-}_{(0,1)}\otimes t\big)^{(d)}v&\\
&=\big(x^{-}_{(1,0)}\otimes t\big)^{(a)}\big(x^{-}_{(0,1)}\otimes t\big)^{(b+c)}\big(x^{-}_{(1,0)}\otimes 1\big)^{(b+2c)}\big(x^{-}_{(0,1)}\otimes t\big)^{(d)}v.
\end{align*}
Now passing all $\big(x_{(0,1)}^-\otimes t\big)$ through $\big(x_{(1,0)}^{-}\otimes 1\big)$ we obtain as above
$$0=\sum_{} c_{i,j}\big(x^{-}_{(1,0)}\otimes t\big)^{(a)}\big(x^{-}_{(1,0)}\otimes 1\big)^{(b+2c-i-2j)}\big(x^{-}_{(1,1)}\otimes t\big)^{(i)}\big(x^{-}_{(2,1)}\otimes t\big)^{(j)}\big(x^{-}_{(0,1)}\otimes t\big)^{(b+c+d-i-j)}v,$$
where the sum runs again over all tuples of non-negative integers $(i,j)$ such that $i+j\leq b+c$ and $i+2j\leq b+2c$ and each $c_{i,j}$ is non-zero.
Now we want to find the maximal monomial in the above sum. Since we want to minimize the power of $\big(x_{(0,1)}^{-}\otimes t\big)$ we should have $i+j=b+c$ in the maximal monomial. If $i<b$ in the maximal monomial, we would get $j>c$ and hence
$$i+2j=b+c+j>b+2c$$
which is a contradiction. Hence $i=b$ by the definition of the total order. This implies $j=c$ and the maximal monomial is given as desired.
\end{proof}
\end{prop}
Again Proposition~\ref{pr2} finishes all parts of Proposition~\ref{wichtp} except the statement on the cardinality of $\mathcal{S}_{\lambda,\mu}^{C}$.
\subsection{}Again there are several well-known combinatorial models for the sum of Littlewood-Richardson coefficients; none of the being compatible with the fusion product grading. We consider here the set $\mathcal{T}^{C}_{\lambda,\mu}$ consisting of $(a,b,c,d) \in \mathbb{Z}_{+}^4$ satisfying the following inequalities: 
\begin{gather*}
		a \le m_1, \qquad c \le m_2 , \qquad d \le n_2, \qquad b \le n_1, \\ 
		c+b-a \le m_2, \qquad d + b -a \le m_2 , \qquad 
		a + 2(c-d) \le n_1, \qquad b+ 2(c-d) \le n_1.
	\end{gather*}	
Then it has been proved in \cite[Section 3]{BZ88} that 
$$|\mathcal{T}^{C}_{\lambda,\mu}|=\sum_{\nu\in P^+}[V(\lambda)\otimes V(\mu): V(\nu)].$$
We define an injective map 
$$\mathcal{T}^{C}_{\lambda-\varpi_1,\mu-\varpi_1}\rightarrow \mathcal{T}^{C}_{\lambda,\mu},\ \ (a,b,c,d)\mapsto (a+1,b+1,c,d),$$
whose image consists exactly of those points in $\mathcal{T}^{C}_{\lambda,\mu}$ where $a,b>0$. Hence 
\begin{equation}\label{typc1}
|\mathcal{T}^{C}_{\lambda,\mu}|=|\mathcal{T}^{C}_{\lambda-\varpi_1,\mu-\varpi_1}|+|^{1}\mathcal{T}^{C}_{\lambda,\mu}|,
\end{equation}
where $^{1}\mathcal{T}^{C}_{\lambda,\mu}$ is the union of the following two polytopes: 
the one consisting of $(0,b,c,d) \in \mathbb{Z}_{+}^4$ satisfying 
\begin{align*}
	 c+b \le m_2 , \qquad   d \le n_2, \qquad  b \le n_1,\qquad   d + b  \le m_2 , \qquad  b+ 2(c-d) \le n_1 
\end{align*}
and that consisting of $(a,0,c,d) \in \mathbb{Z}_+^4$ satisfying 
\begin{align*}
	a \ge 1,\qquad   a \le m_1, \qquad c \le m_2 , \qquad d \le n_2, \qquad  d  -a \le m_2 ,\qquad  a + 2(c-d) \le n_1.
\end{align*}
Using the injective map 
\begin{align*}
	^{1}\mathcal{T}^{C}_{\lambda-\varpi_2,\mu-\varpi_2} \rightarrow {}^{1}\mathcal{T}^{C}_{\lambda,\mu}, \ \ (a,b,c,d) \mapsto (a,b,c+1,d+1)
\end{align*}
we further obtain 
\begin{equation}\label{typc1:t2}
	|^{1}\mathcal{T}^{C}_{\lambda,\mu}|=|^{1}\mathcal{T}^{C}_{\lambda-\varpi_2,\mu-\varpi_2}|+|^{2}\mathcal{T}^{C}_{\lambda,\mu}|,
\end{equation}
where $^{2}\mathcal{T}^{C}_{\lambda, \mu}\subseteq \mathbb{Z}_+^4$ is the disjoint union of the following four polytopes:
\begin{align*}
		&_{1}^{2}\mathcal{T}^{C}_{\lambda, \mu}:= \left\{ (0,b,0,d): 
		\substack{\displaystyle  d  \le n_2, \, b \le n_1\\ \\ \displaystyle  d + b  \le m_2  }\right\}, 
		&&
		_{2}^{2}\mathcal{T}^{C}_{\lambda, \mu}:=\left\{ (a,0,0,d): 
		\substack{\displaystyle  1 \le a \le m_1, \,   d \le n_2\\ \\ \displaystyle  d  -a \le m_2 , \, a -2d \le n_1    }\right\}, 	
		\\ 
		&_{3}^{2}\mathcal{T}^{C}_{\lambda, \mu}:= \left\{ (0,b,c,0): 
		\substack{\displaystyle  c \ge 1, \, c +b\le m_2 \\ \\ \displaystyle  b+ 2c \le n_1    }\right\}, 
		&&
		_{4}^{2}\mathcal{T}^{C}_{\lambda, \mu}:=\left\{  (a,0,c,0) : 
		\substack{\displaystyle 1 \le a \le m_1, \, 1 \le c \le m_2 \\ \\ \displaystyle  a + 2c \le n_1    }\right\}. 
\end{align*}

We summarize the combinatorics of $\mathcal{T}^{G}_{\lambda,\mu}$ in the following lemma.
\begin{lem}\label{334c} Let $\lambda,\mu\in P^+$. Then 
$$|\mathcal{T}^{C}_{\lambda,\mu}|=|\mathcal{T}^{C}_{\lambda-\varpi_1,\mu-\varpi_1}|+|^{1}\mathcal{T}^{C}_{\lambda-\varpi_2,\mu-\varpi_2}|+|^{2}\mathcal{T}^{C}_{\lambda,\mu}|.$$
Moreover, if $\min\{m_2, n_2\} > 0$, then 
$$|^{2}\mathcal{T}^{C}_{\lambda, \mu}| =|^{2}\mathcal{T}^{C}_{\lambda-\varpi_2, \mu-\varpi_2}| + \min\left\{2(m_1 + m_2), 2(n_1 + n_2), m_1 + n_1\right\} +1$$
and if $\min\{m_2, n_2\} =0$ and $\min\{m_1, n_1\} > 0$ we have
$$ |^{2}\mathcal{T}^{C}_{\lambda-\varpi_1, \mu-\varpi_1}| + \min\{n_1, m_2\} + \min\{m_1, n_2\} +1.$$
\proof The first part follows from \eqref{typc1} and  \eqref{typc1:t2}.  Assume that $\min\{m_2, n_2\} > 0$. Analyzing the polytopes defining 
	$^{2}\mathcal{T}^{C}_{\lambda, \mu}$, one can calculate the following differences:
	\begin{align*}
		A_1 := |{}_{1}^{2}\mathcal{T}^{C}_{\lambda, \mu}| - |{}_1^{2}\mathcal{T}^{C}_{\lambda-\varpi_2, \mu-\varpi_2}|	
		& = 1+ \min\{n_1, m_2\},\\ 
		A_2:= |{}_{2}^{2}\mathcal{T}^{C}_{\lambda, \mu}| - |{}_2^{2}\mathcal{T}^{C}_{\lambda-\varpi_2, \mu-\varpi_2}|
		 &=(\min\{m_1, n_1 + 2n_2\} - (n_2 - m_2)_+)_+ + \min\{m_1, n_2 - m_2\}_+,\\
		A_3 := |{}_{3}^{2}\mathcal{T}^{C}_{\lambda, \mu}| - |{}_3^{2}\mathcal{T}^{C}_{\lambda-\varpi_2, \mu-\varpi_2}|	
		& = \min\{n_1 - m_2, m_2\}_+,\\
		A_4 := |{}_{4}^{2}\mathcal{T}^{C}_{\lambda, \mu}| - |{}_4^{2}\mathcal{T}^{C}_{\lambda-\varpi_2, \mu-\varpi_2}|	
		& = \min\{m_1, n_1 - 2m_2\}_+.
	\end{align*}
	This leads to 
	\begin{align*}
		A_1 + A_3 + A_4 & = 
		\begin{cases}
			n_1, &\text{if } n_1 < 2m_2, \\ 
			\min\{n_1, m_1 + 2m_2\}, & \text{if }n_1 \ge 2m_2.
		\end{cases}\\ 
		A_2 & = \begin{cases}
			m_1, &\text{if } m_1 < n_2-m_2, \\ 
			\min\{m_1, n_1 + 2n_2\}, & \text{if }m_1 \ge n_2-m_2.
		\end{cases}
	\end{align*}
	This implies the second statement:
	\begin{align*}
		|^{2}\mathcal{T}^{C}_{\lambda, \mu}| - | ^{2}\mathcal{T}^{C}_{\lambda-\varpi_2, \mu-\varpi_2}|  = 
		A_1 + A_2 + A_3 + A_4 
		= \min\{2(m_1 + m_2), 2(n_2 + n_2),m_1+n_1\}+1.
	\end{align*}
Now suppose $\min\{m_2,n_2\}=0$ and $\min\{m_1,n_1\}>0$. If $n_2 = 0$, then $^{2}\mathcal{T}^{C}_{\lambda, \mu}$ consists of the union of the following polytopes
	\begin{align*}
	\begin{split}
		\left\{ (0,b,0,0): 
		\substack{\displaystyle b \le \min\{n_1,m_2\}}\right\} ,
		\qquad & \qquad \ \quad 
		\left\{ (0,b,c,0): 
		\substack{\displaystyle  c \ge 1, \, c +b\le m_2 \\ \\ \displaystyle  b+ 2c \le n_1    }\right\} ,
		\\
		\left\{  (a,0,c,0) : 
		\substack{\displaystyle 1 \le a \le m_1, \, 1 \le c \le m_2 \\ \\ \displaystyle  a + 2c \le n_1    }\right\} ,
		&\qquad 
		\left\{ (a,0,0,0): 
		\substack{\displaystyle  1 \le a \le \min\{n_1,m_1\}}\right\},
	\end{split}
\end{align*}
which can be unified to the union of the following two polytopes:
\begin{align*}
		\left\{ (0,b,c,0): 
		\substack{\displaystyle  c +b\le m_2 \\ \\ \displaystyle  b+ 2c \le n_1    }\right\} ,
		&&
		\left\{  (a,0,c,0) : 
		\substack{\displaystyle 1 \le a \le m_1, \, c \le m_2 \\ \\ \displaystyle  a + 2c \le n_1    }\right\}.
\end{align*}
Thus $|^{2}\mathcal{T}^{C}_{\lambda, \mu}| - |{}^{2}\mathcal{T}^{C}_{\lambda-\varpi_1, \mu-\varpi_1}|$ is given by the cardinality of the disjoint union
\begin{align*}
		\left\{ (0,b,c,0): 
		\substack{\displaystyle  c +b\le m_2 \\ \\ \displaystyle  b+ 2c = n_1    }\right\} ,
		&& 
		\left\{  (1,0,c,0) : 
		\substack{\displaystyle  c \le m_2 \\ \\ \displaystyle  2c \le n_1-1   }\right\} ,
\end{align*}

which is equal to 
	\begin{align*}
		\left(\left\lfloor\frac{n_1}{2}\right\rfloor  - (n_1 - m_2)_+ +1\right)_++ \min\left\{ \left\lfloor\frac{n_1-1}{2}\right\rfloor , m_2\right\} +1 
		= \min\{n_1, m_2\} +1.
	\end{align*}
The first summand describes the cardinality of the first polytope and the second summand the cardinality of the second one. Finally if $m_2 = 0$, then the cardinality of $^{2}\mathcal{T}^{C}_{\lambda, \mu}$ is given by $|_{2}^{2}\mathcal{T}^{C}_{\lambda, \mu}|+1$, where $_{2}^{2}\mathcal{T}^{C}_{\lambda, \mu}$ simplifies to
	$$
		_{2}^{2}\mathcal{T}^{C}_{\lambda, \mu}=\left\{ (a,0,0,d): 
		\substack{\displaystyle  1 \le a \le m_1, \,   d \le n_2\\ \\ \displaystyle  d  \le a , \, a -2d \le n_1    }\right\}.$$
Now it is straightforward to see that$ |^{2}\mathcal{T}^{C}_{\lambda, \mu}| - | ^{2}\mathcal{T}^{C}_{\lambda-\varpi_1, \mu-\varpi_1}|$ is the same as the cardinality of those elements in $_{2}^{2}\mathcal{T}^{C}_{\lambda, \mu}$ where $a=1$ or ($a>1$ and $d=a$). 
Hence the difference is given by $\min\{n_2, m_1\} + 1$.
\endproof
\end{lem}
\subsection{}Now let us find analoges of the previous statements for $\mathcal{S}^{C}_{\lambda,\mu}$. We can define an injective map 
$$\mathcal{S}^{C}_{\lambda-\varpi_1,\mu-\varpi_1}\rightarrow \mathcal{S}^{C}_{\lambda,\mu},\ \ (a,b,c,d)\mapsto (a+1,b,c,d),$$
whose image consists exactly of those points in $\mathcal{S}^{C}_{\lambda,\mu}$ where $a>0$. Thus 
\begin{equation}\label{typc1:s1}
|\mathcal{S}^{C}_{\lambda,\mu}|=|\mathcal{S}^{C}_{\lambda-\varpi_1,\mu-\varpi_1}|+|^{1}\mathcal{S}^{C}_{\lambda,\mu}|,
\end{equation} 
where $^{1}\mathcal{S}^{C}_{\lambda,\mu}$ consists of $(b,c,d) \in \mathbb{Z}_{+}^3$ satisfying 
\begin{gather*}
	d \le \min\{m_2, n_2\} , \qquad 
	b+d \le \min\{m_1 + m_2, n_1 + n_2\}, \qquad b+c \le \min\{m_1 + m_2, n_1 + n_2\} \\
	  b \le m_1 + n_1, \qquad 2d + b \le m_2 + n_2 , \qquad b +2c - 2d \le m_1 + n_1. 
\end{gather*}
Similar as before, we can define an injective map 
\begin{align*}
	^{1}\mathcal{S}^{C}_{\lambda-\varpi_2,\mu-\varpi_2} \rightarrow {}^{1}\mathcal{S}^{C}_{\lambda,\mu}, \ \ (b,c,d) \mapsto (b,c+1,d+1),
\end{align*}
which leads to 
\begin{equation}\label{typc1:s2}
	|^{1}\mathcal{S}^{C}_{\lambda,\mu}|=|^{1}\mathcal{S}^{C}_{\lambda-\varpi_2,\mu-\varpi_2}|+|^{2}\mathcal{S}^{C}_{\lambda,\mu}|,
\end{equation}
where $^{2}\mathcal{S}^{C}_{\lambda,\mu} \subseteq \mathbb{Z}_+^{3}$ is union of the following two sets: the one 
consisting of $(b,0,d)$ such that 
\begin{align*}
	d \le \min\{m_2, n_2\} , \quad b+d \le \min\{m_1 + m_2, n_1 + n_2\}, \quad b \le m_1 + n_1, \quad 2d + b \le m_2 + n_2
\end{align*}
and the one consisting of $(b,c,0)$ satisfying 
\begin{align*}
	c \ge 1,\quad  b+c \le \min\{m_1 + m_2, n_1 + n_2\} , \quad b \le m_2 + n_2 , \quad b +2c  \le m_1 + n_1.
\end{align*}
\begin{lem}\label{334c2} Let $\lambda,\mu\in P^+$. Then 
$$|\mathcal{S}^{C}_{\lambda,\mu}|=|\mathcal{S}^{C}_{\lambda-\varpi_1,\mu-\varpi_1}|+|^{1}\mathcal{S}^{C}_{\lambda-\varpi_2,\mu-\varpi_2}|+|^{2}\mathcal{S}^{C}_{\lambda,\mu}|.$$
Moreover, if $\min\{m_2, n_2\} > 0$, then 
$$|^{2}\mathcal{S}^{C}_{\lambda, \mu}| =|^{2}\mathcal{S}^{C}_{\lambda-\varpi_2, \mu-\varpi_2}| + \min\left\{2(m_1 + m_2), 2(n_1 + n_2), m_1 + n_1\right\} +1$$
and if $\min\{m_2, n_2\} =0$ and $\min\{m_1, n_1\} > 0$ we have
$$ |^{2}\mathcal{S}^{C}_{\lambda-\varpi_1, \mu-\varpi_1}| + \min\{n_1, m_2\} + \min\{m_1, n_2\} +1.$$
\proof The first part follows from \eqref{typc1:s1} and  \eqref{typc1:s2}. Assume first that $\min\{m_2, n_2\} > 0$. We have an injective map 
	\begin{align*}
		^{2}\mathcal{S}^{C}_{\lambda-\varpi_2, \mu-\varpi_2} 
		\rightarrow 
		{}^{2}\mathcal{S}^{C}_{\lambda, \mu}, 
		(b,c,d) 
		\mapsto 
		\begin{cases}
			(b,0,d+1), & \text{if } c = 0,\\
			(b+2,c-1,0) & \text{otherwise},
		\end{cases}
	\end{align*}
whose image does not contain the elements in ${}^{2}\mathcal{S}^{C}_{\lambda, \mu}$ of the form $(b,c,0)$ with $b\leq 1$ and every other point is in the image. Hence we obtain that
	\begin{align*}
		|{}^{2}\mathcal{S}^{C}_{\lambda, \mu}|= |{}^{2}\mathcal{S}^{C}_{\lambda-\varpi_2, \mu-\varpi_2}|+|{}^{3}\mathcal{S}^{C}_{\lambda, \mu}|,
	\end{align*}
	where ${}^{3}\mathcal{S}^{C}_{\lambda, \mu}$ consists of all points $(b,c)\in\mathbb{Z}^2_+$ satisfying
	$$b\leq 1,\quad b+c \le \min\{m_1 + m_2, n_1 + n_2\} , \quad b +2c  \le m_1 + n_1.$$
Since
	\begin{align*}
		|{}^{3}\mathcal{S}^{C}_{\lambda, \mu}|=\min\{2(m_1 + m_2), 2(n_1 + n_2), m_1 + n_1 \} + 1,
	\end{align*}
	the claim follows in this case.
If $\min\{m_2, n_2\} = 0$, then ${}^{2}\mathcal{S}^{C}_{\lambda, \mu}$ consists of all $(b,c,0) \in \mathbb{Z}_{+}^3$ satisfying 
	\begin{align*}
		b+c \le \min\{m_1 + m_2, n_1 + n_2\} , \quad b \le m_2 + n_2 , \quad b +2c  \le m_1 + n_1.
	\end{align*}
	With the injective map
	\begin{align*}
		^{2}\mathcal{S}^{C}_{\lambda-\varpi_1, \mu-\varpi_1} \rightarrow {}^{2}\mathcal{S}^{C}_{\lambda, \mu}, \
		(b,c,0) \mapsto (b,c+1,0)
	\end{align*}
	we easily see that $|{}^{2}\mathcal{S}^{C}_{\lambda, \mu}| - |^{2}\mathcal{S}^{C}_{\lambda-\varpi_1, \mu-\varpi_1}| $ equals $\min\{m_1, n_2\} + \min\{n_1, m_2\} +1.$
\endproof

\end{lem}
\subsection{}
So the following lemma finishes the proof of Proposition~\ref{wichtp} in the case of $C_2$.
\begin{lem}\label{}
There exists a bijective map $\mathcal{S}^{C}_{\lambda,\mu} \rightarrow  \mathcal{T}^{C}_{\lambda,\mu}$.
\begin{proof}
In view of Lemma \ref{334c}, Lemma \ref{334c2} and an easy induction argument we only have to prove that $|{}^{2}\mathcal{T}^{C}_{\lambda, \mu}| =|{}^{2}\mathcal{S}^{C}_{\lambda, \mu}|$ if $\min\{m_1, n_1\} = 0=\min\{m_2, n_2\}$. If $n_1, n_2 = 0$ or $m_1, m_2 = 0$, then both sets consist only of one point. If $m_1 = 0 = n_2$, then both sets are defined by the same non redundant inequalities, so the only interesting case is $m_2 = 0 = n_1$ and $n_2 > 0$. Now the statement follows from 
	\begin{align*}
		|\mathcal{S}^{C}_{\lambda,\mu}|- |\mathcal{S}^{C}_{\lambda,\mu-\varpi_2}| =   (\min\{n_2, m_1 - n_2 \}+1)_+ =
		|\mathcal{T}^{C}_{\lambda,\mu}|- |\mathcal{T}^{C}_{\lambda,\mu-\varpi_2}|,
	\end{align*}
	which can be easily seen by arguments similar to those we used before. 
\end{proof}
\end{lem}
\section{The case $G_2$}\label{section6}
In this section we prove Proposition~\ref{wichtp} for the rank two Lie algebra of type $G_2$.
\subsection{}Here we consider the total order on $\mathbb{Z}^6_+$ defined as follows $\mathbf{s}=(a,b,c,d,e,f)\succ \mathbf{s}'=(a',b',c',d',e',f')$ if and only if one of the following conditions hold
$$\circ\ f<f',\quad
\circ\ f=f',\ a<a'  \quad
\circ\  f=f',\ a=a',\ b<b'\quad
\circ\ f=f',\ a=a',\ b=b',\ c<c'\quad
$$
$$
\circ\  f=f',\ a=a',\ b=b',\ c=c',\ d<d'\quad
\circ\ f=f',\ a=a',\ b=b',\ c=c',\ d=d',\ e<e'.$$
\begin{lem}\label{eachsumm}
The following relations hold in $\mathcal{F}_{\lambda,\mu}$ provided that $\min\{m_2,n_2\}=0$:
$$\big(x^{-}_{(3,1)}\otimes t\big)^{(s)}\big(x^{-}_{(1,0)}\otimes 1\big)^{(r)}v=0,\ \ \ \forall r,s\in\mathbb{Z}_+:  r+s>m_1+n_1$$
$$\big(x^{-}_{(1,0)}\otimes 1\big)^{(r)}\big(x^{-}_{(1,1)}\otimes t\big)^{(s)}v=0,\ \ \ \forall r,s\in\mathbb{Z}_+:  r-s>m_1+n_1.$$
\begin{proof}
We prove the lemma by induction on $s$, where the case $s=0$ is clear. Now if $s>0$ we use \cite[Lemma 4.1]{SaCh} and obtain
\begin{align}
\notag0&=\big(x^{-}_{(1,0)}\otimes 1\big)^{(r+3s)}\big(x^{-}_{(0,1)}\otimes t\big)^{(s)}v&\\&\label{eachsumm1}
=\sum_{}\big(x^{-}_{(0,1)}\otimes t\big)^{(s-i-j-k)}\big(x^{-}_{(1,1)}\otimes t\big)^{(i)}\big(x^{-}_{(2,1)}\otimes t\big)^{(j)}\big(x^{-}_{(3,1)}\otimes t\big)^{(k)}\big(x^{-}_{(1,0)}\otimes 1\big)^{(\ell)}v,
\end{align}
where the sum runs over all $(i,j,k,\ell)$ of non-negative integers satisfing $i+j+k\leq s$ and $i+2j+3k+\ell=r+3s$. Since 
$$\ell+k=(r+3s-i-2j-2k)\geq r+s>m_1+n_1,$$ 
we obtain by induction that each summand for $k<s$ in \eqref{eachsumm1} is zero. Now the claim follows by taking $k=s$ in \eqref{eachsumm1}. In order to see the second relation we get by using \cite[Lemma 4.1]{SaCh}:
$$\big(x^{-}_{(1,0)}\otimes 1\big)^{(r)}\big(x^{-}_{(1,1)}\otimes t\big)^{(s)}v=\sum_{\substack{i+j\leq s,\\ i+2j\leq r}}\big(x^{-}_{(1,1)}\otimes t\big)^{(s-i-j)}\big(x^{-}_{(2,1)}\otimes t\big)^{(i)}\big(x^{-}_{(3,1)}\otimes t\big)^{(j)}\big(x^{-}_{(1,0)}\otimes 1\big)^{(r-i-2j)}=0.$$
Now the claim follows immediately from the first relation.
\end{proof}
\end{lem}
\begin{prop}\label{pr3}
Assume that $\mathbf{s}\in \mathbb{Z}_+^6$ is such that $\mathbf{s}\notin \mathcal{S}^G_{\lambda,\mu}$ Then we have
$$\left(X_{\mathbf{s}}+\sum_{\mathbf{s}'\prec\mathbf{s} } c_{\mathbf{s}'} X_{\mathbf{s}'}\right)v=0.$$
for suitable elements $c_{\mathbf{s}'}\in\mathbf{U}(\mathfrak{n}^-)$.
\begin{proof}
If $\mathbf{s}\in \mathbb{Z}_+^6$ violates one of the equalities in \eqref{6}, then the statement is clear by the defining equalities of $\mathcal{F}_{\lambda,\mu}$ or by the following argument. If $b>m_2+n_2$, then 
\begin{align*}
	0=\big(x_{(1,0)}^-\otimes t\big)^{(a+b)}&\big(x_{(0,1)}^-\otimes 1\big)^{(b)}v=\sum^b_{i=0} \pm\big(x_{(0,1)}^{-}\otimes 1\big)^{(b-i)}(x^{-}_{(1,1)}\otimes t\big)^{(i)}\big(x^{-}_{(1,0)}\otimes t\big)^{(a+b-i)}v,
\end{align*}
where the sum is made by passing successively $\big(x_{(1,0)}^-\otimes t\big)$ through $\big(x_{(0,1)}^-\otimes 1\big)$ (one can check this directly or use \cite[Lemma 4.1]{SaCh}). 
By the definition of the total order, the maximal monomial in the above sum arises when $i=b$ and in each other monomial the power of $\big(x_{(0,1)}^-\otimes t\big)$ is strictly greater than $a$. Now the rest follows by multiplying the above equation by $\big(x^{-}_{(2,1)}\otimes t\big)^{(c)}\big(x^{-}_{(3,1)}\otimes t\big)^{(d)}\big(x^{-}_{(3,2)}\otimes t\big)^{(e)}$ since passing this product through $\big(x_{(0,1)}^-\otimes 1\big)$ has no affect on the power of $\big(x_{(0,1)}^-\otimes t\big)$. To see the last inequality in \eqref{6}, we consider in the case when $b+e-a>m_2+n_2$ the product
$$0=\big(x^-_{(0,1)}\otimes 1\big)^{(b+e)}\big(x^-_{(1,0)}\otimes t\big)^{(a)}v.$$
Then we get (again by direct calculation or using twice \cite[Lemma 4.1]{SaCh}) 
\begin{align*}
	0&=\big(x^-_{(1,0)}\otimes t\big)^{(b)}\big(x^-_{(3,1)}\otimes t\big)^{(d+e)}\big(x^-_{(0,1)}\otimes 1\big)^{(b+e)}\big(x^-_{(1,0)}\otimes t\big)^{(a)}v&\\
	&=\sum_{i=0}^{e+\min\{b,d\}}\pm \big(x^-_{(1,0)}\otimes t\big)^{(b)}\big(x^-_{(0,1)}\otimes 1\big)^{(b+e-i)}\big(x^-_{(3,2)}\otimes t\big)^{(i)}
	\big(x^-_{(3,1)}\otimes t\big)^{(d+e-i)}\big(x^-_{(1,0)}\otimes t\big)^{(a)}v&\\
	&=\sum_{i=0}^{e+\min\{b,d\}}\sum_{j=0}^{b+\min\{0,e-i\}}c_j \big(x^-_{(0,1)}\otimes 1\big)^{(b+e-i-j)}\big(x^-_{(1,1)}\otimes 1\big)^{(j)} &\\&
	 \hspace{4cm} \times \big(x^-_{(3,2)}\otimes t\big)^{(i)}\big(x^-_{(3,1)}\otimes t\big)^{(d+e-i)}\big(x^-_{(1,0)}\otimes t\big)^{(a+b-j)}v,\ \ \ c_j\neq 0.
\end{align*}
The maximal monomial arises when $i=e$ and $j=b$ and the claim follows by multiplying the above sum with $\big(x^-_{(2,1)}\otimes t\big)^{(c)}$ (since $x^{-}_{(2,1)}$ and $x^{-}_{(0,1)}$ commute). So suppose in the rest of the proof that $\mathbf{s}$ satisfies the equations in \eqref{6}. In the next step we suppose that the first inequality in \eqref{7} is violated (the proof for the second inequality in \eqref{7} and the inequalities in \eqref{71} is similar and will be omitted). 
We get 
$$0=\big(x^{+}_{(1,0)}\otimes 1\big)^{(c)}\big(x^{+}_{(2,1)}\otimes 1\big)^{(a)}\big(x^{-}_{(3,1)}\otimes t\big)^{(a+c+d)}v.$$ 
Since $\big(x_{(1,1)}^+\otimes t\big)\big(x^{+}_{(2,1)}\otimes 1\big)^{(r)}\big(x^{-}_{(3,1)}\otimes t\big)^{(s)}v=0$ for all $r,s\in\mathbb{Z}_+$ we get
\begin{align*}0=\big(x^{+}_{(2,1)}\otimes 1\big)^{(a)}\big(x^{-}_{(3,1)}\otimes t\big)^{(a+c+d)}v&=\big(x^{+}_{(2,1)}\otimes 1\big)^{(a-1)}\big(x^{-}_{(1,0)}\otimes t\big)\big(x^{-}_{(3,1)}\otimes t\big)^{(a+c+d-1)}v&\\&=\big(x^{-}_{(1,0)}\otimes t\big)\big(x^{+}_{(2,1)}\otimes 1\big)^{(a-1)}\big(x^{-}_{(3,1)}\otimes t\big)^{(a+c+d-1)}v&\\&=\big(x^{-}_{(1,0)}\otimes t\big)^{(a)}\big(x^{-}_{(3,1)}\otimes t\big)^{(c+d)}v.\end{align*}
Hence
\begin{align*}0=\big(x^{+}_{(1,0)}\otimes  1\big)^{(c)}&\big(x^{-}_{(1,0)}\otimes t\big)^{(a)}\big(x^{-}_{(3,1)}\otimes t\big)^{(c+d)}v=\big(x^{-}_{(1,0)}\otimes t\big)^{(a)}\big(x^{+}_{(1,0)}\otimes 1\big)^{(c)}\big(x^{-}_{(3,1)}\otimes t\big)^{(c+d)}v&\\&=\big(x^{-}_{(1,0)}\otimes t\big)^{(a)}\left(\big(x^{-}_{(2,1)}\otimes t\big)^{(c)}\big(x^{-}_{(3,1)}\otimes t\big)^{(d)}v+\sum_i X_i \right)v, \end{align*}
where each monomial $X_i$ in the above sum has one of the elements $\big(x^{-}_{(0,1)}\otimes t\big)$
or $\big(x^{-}_{(1,1)}\otimes t\big)$ as a factor. So the claim follows by multiplying the above equation by $\big(x^{-}_{(1,1)}\otimes t\big)^{(b)}\big(x^{-}_{(3,2)}\otimes t\big)^{(e)}$ and the definition of the order. 
We finally consider the case when one of the inequalities in \eqref{8} are violated, starting with $2(a+c)+3d-b>m_1+n_1$. Up to sign, we get from Lemma~\ref{eachsumm}
\begin{align}\notag 0&=\big(x^{+}_{(1,0)}\otimes t\big)^{(a)}\big(x^{-}_{(1,0)}\otimes 1\big)^{(2(a+c)+3d)}\big(x^{-}_{(1,1)}\otimes t\big)^{(b)}v&\\& \label{eachsumm11}=\big(x^{-}_{(1,0)}\otimes t\big)^{(a)}\big(x^{-}_{(1,0)}\otimes 1\big)^{(2c+3d)}\big(x^{-}_{(1,1)}\otimes t\big)^{(b)}v.\end{align}
Now multiplying \eqref{eachsumm11} with $\big(x^-_{(0,1)}\otimes t\big)^{(c+d)}$ gives (using once more \cite[Lemma 4.1]{SaCh})
\begin{align*}0&=\big(x^-_{(1,0)}\otimes t\big)^{(a)}\big(x^-_{(0,1)}\otimes t\big)^{(c+d)}\big(x^{-}_{(1,0)}\otimes 1\big)^{(2c+3d)}\big(x^{-}_{(1,1)}\otimes t\big)^{(b)}v&\\&=\sum_{}c_{i,j,k}\big(x^-_{(1,0)}\otimes 1\big)^{(2c+3d-i-2j-3k)}\big(x^-_{(1,0)}\otimes t\big)^{(a)}&\\&\hspace{5cm}\times \big(x^-_{(1,1)}\otimes t\big)^{(b+i)}\big(x^-_{(2,1)}\otimes t\big)^{(j)}\big(x^-_{(3,1)}\otimes t\big)^{(k)}v,\ \ c_{i,j,k}\neq 0,\end{align*}
where the sum runs over triples $(i,j,k)$ satisfying $i+j+k=c+d$ and $i+2j+3k\leq 2c+3d$. 
By the definition of the order we want to minimalize $i$, so choose $i=0$ in the above sum. If $j<c$ we would get $k>d$ and hence 
$$2j+3k=2(j+k)+k=2c+2d+k>2c+3d,$$which is a contradiction. Hence the maximal monomial appears if we choose $i=0$, $j=c$ and $k=d$. 
In the last case we assume $2(a+c)+b+d>m_1+n_1$. From Lemma~\ref{eachsumm} we know that up to sign 
$$0=\big(x^+_{(1,0)}\otimes t\big)^{(a)}\big(x^-_{(1,0)}\otimes 1\big)^{(2(a+c)+b)}\big(x^-_{(3,1)}\otimes t\big)^{(d)}v=\big(x^-_{(1,0)}\otimes t\big)^{(a)}\big(x^-_{(1,0)}\otimes 1\big)^{(2c+b)}\big(x^-_{(3,1)}\otimes t\big)^{(d)}v.$$
Multiplying the above equation with $\big(x^-_{(0,1)}\otimes t\big)^{(b+c)}$ we get exactly as above 
\begin{align*}0&=\big(x^-_{(1,0)}\otimes t\big)^{(a)}\big(x^-_{(0,1)}\otimes t\big)^{b+c}\big(x^-_{(1,0)}\otimes 1\big)^{(2c+b)}\big(x^-_{(3,1)}
\otimes t\big)^{d}v&\\&=\sum_{}c_{i,j,k}\big(x^-_{(1,0)}\otimes 1\big)^{(2c+b-i-2j-3k)}\big(x^-_{(1,0)}\otimes t\big)^{(a)} &\\& \hspace{4cm} \times \big(x^-_{(1,1)}\otimes t\big)^{(i)}\big(x^-_{(2,1)}\otimes t\big)^{(j)}\big(x^-_{(3,1)}\otimes t\big)^{(d+k)}v,\ \ c_{i,j,k}\neq 0,\end{align*}
where the sum runs over triples $(i,j,k)$ satisfying $i+j+k=b+c$ and $i+2j+3k\leq 2c+b$. 
The conditions on $(i,j,k)$ imply $c \ge j+2k = b+c-i+k$ hence $k \le i-b$ and $i\ge b$. The
maximal element is therefore given by choosing $i = b, k = 0$ and $j = c$.
This finishes the proof of the proposition.
\end{proof}
\end{prop}
Again it remains to investigate the cardinality of $\mathcal{S}_{\lambda,\mu}^{G}$.
\subsection{} In the rest of this section we assume for simplicity that $n_2=0$ and consider the set $\mathcal{T}^{G}_{\lambda,\mu}$ consisting of all elements $(a,b,c,d,e,f) \in \mathbb{Z}_{+}^6$ satisfying the following inequalities
\begin{gather*}
		a+b+c+d+e+f  \le n_1, \qquad c \le 1, \qquad b+e-d\le m_2, \label{eq:lit1} \\
		f\leq m_1,\ \qquad e\leq m_2,\qquad  a-2b+2d-e+f  \le m_1,  \label{eq:lit2}
		\\
		c+f+2d-e\leq m_1 \label{eq:lit3}
	\end{gather*}	
\begin{lem}\label{Li1}
We have an equality
$$|\mathcal{T}^{G}_{\lambda,\mu}|=\sum_{\nu\in P^+}[V(\lambda)\otimes V(\mu): V(\nu)].$$
\begin{proof}We consider the following admissible sixtuples of rows
$$\begin{ytableau}
1  \\
1\\
1  \\
1\\
1  \\
1
\end{ytableau} \quad \quad \quad 
\begin{ytableau}
2  \\
2\\
2  \\
2\\
2  \\
2
\end{ytableau} \quad \quad \quad 
\begin{ytableau}
3  \\
3\\
3  \\
3\\
3  \\
3
\end{ytableau} 
 \quad \quad \quad 
\begin{ytableau}
3  \\
3\\
3  \\
4\\
4  \\
4
\end{ytableau} 
 \quad \quad \quad 
\begin{ytableau}
4  \\
4\\
4  \\
4\\
4  \\
4
\end{ytableau} 
 \quad \quad \quad 
\begin{ytableau}
5  \\
5\\
5  \\
5\\
5  \\
5
\end{ytableau} 
 \quad \quad \quad 
\begin{ytableau}
6  \\
6\\
6  \\
6\\
6  \\
6
\end{ytableau} 
$$
and call them $Q_1,Q_2,Q_3,Q_{3,4},Q_4,Q_5$ and $Q_6$ respectively. Let $T$ be a column tableau of length $6n_1$ such that the entries in the boxes are not decreasing from the top to the bottom and for all $\ell=1,\dots,n_1$ the subtableau consisting of the $(6\ell-5)$th up to the $6\ell$th row is an admissible sixtuple of rows. Let $T_i$ be the subtablau obtained from $T$ by considering the last $i$ rows. Define
$$c_{i,j}=\text{number of $j$'s in the tableau $T_i$}$$
and let 
$$d^{i}_1=c_{i,1}+2c_{i,3}+c_{i,5}-c_{i,2}-c_{i,6}-2c_{i,4}+6m_1,\ \ \  d^{i}_2=c_{i,2}+c_{i,4}-c_{i,5}-c_{i,3}+6m_2.$$
If $d^{i}_1,d^{i}_1\geq 0$ for all $i=1,\dots,6n_1$, we call $T$ to be a standard dominant tableau. From \cite[Theorem 3.7, Section 3.8]{Li90} we get that the number of standard dominant tableaux coincides with the cardinality of 
$$\sum_{\nu\in P^+}[V(\lambda)\otimes V(\mu): V(\nu)].$$
In order to finish the proof of the lemma, let $T$ be a standard dominant tableau and note that $T$ must be a single column starting with a certain number of $Q_1's$ followed by a certain number of $Q_2's$ and so on.
We denote by $y_1$ the number of $Q_1's$ in $T$, by $y_2$ the number of $Q_2's$ in $T$ and so on. We will now investigate step by step the condition being a standard dominant tableau. This will lead to certain inequalities among the non-negative integers $y_1,y_2,y_3,y_{3,4},y_4,y_5,y_6$. First of all, we must have $y_{3,4}\leq 1$ since otherwise the entries in the boxes are decreasing. Since the column length is exactly $6n_1$ we also get the inequality 
$$y_2+y_3+y_{3,4}+y_4+y_5+y_6\leq n_1.$$
Moreover, choosing 
$$i_1=6y_6,\ \ i_2=6(y_6+y_5),\ \ i_3=6(y_6+y_5+y_4+y_3),$$
$$i_4=6(y_3+y_{3,4}+y_4+y_5+y_6),\ \ i_5=6(y_{3,4}+y_4+y_5+y_6),$$
gives the inequalities
$$d_1^{i_1}=6(m_1-y_6)\geq 0,\ \ d_2^{i_2}=6(m_2-y_5)\geq 0,\ \ d_2^{i_3}=6(m_2-y_5-y_3+y_4)\geq 0,$$
$$d_1^{i_4}=6(m_1-y_{3,4}-2y_4-y_6+y_5)\geq 0,\ \ d_1^{i_5}=6(m_1-2y_4-y_2-y_6+2y_3+y_5)\geq 0.$$
It is straightforward to check that the above conditions are also enough to guarantee that $T$ is a standard dominant tableau. Hence we obtain a bijection between the set of standard dominant tableau and $\mathcal{T}^{G}_{\lambda,\mu}$ which maps an element $T$ to the tuple $(y_2,y_3,y_{3,4},y_4,y_5,y_6)$.
\end{proof}
\end{lem}	
\subsection{}We have an injective map
$$\mathcal{T}^{G}_{\lambda-\varpi_1,\mu-\varpi_1}\rightarrow \mathcal{T}^{G}_{\lambda,\mu},\ (a,b,c,d,e,f)\mapsto (a,b,c,d,e,f+1),$$
whose image consists exactly of those points in $\mathcal{T}^{G}_{\lambda,\mu}$ whose last entry is strictly greater than zero. Hence we obtain 
\begin{equation}\label{recursion1}|\mathcal{T}^{G}_{\lambda,\mu}|=|\mathcal{T}^{G}_{\lambda-\varpi_1,\mu-\varpi_1}|+|^{1}\mathcal{T}^{G}_{\lambda,\mu}|,\end{equation}
where $^{1}\mathcal{T}^{G}_{\lambda,\mu}$ consists of tuples $(a,b,c,d,e) \in \mathbb{Z}_{+}^5$ satisfying the following inequalities 
\begin{gather*}
	\label{frt11}	a+b+c+d+e  \le n_1, \qquad c \le 1, \qquad b+e-d\leq m_2, \\
	\label{frt22}	e\leq m_2,\quad a-2b+2d-e\leq m_1,\qquad c+2d-e \le m_1.   
	\end{gather*}	
Again we can define an injective map 
$$^{1}\mathcal{T}^{G}_{\lambda+\varpi_1-\varpi_2,\mu-\varpi_1}\rightarrow \ ^{1}\mathcal{T}^{G}_{\lambda,\mu},\ (a,b,c,d,e)\mapsto (a,b,c,d,e+1),$$
whose image consists exactly of those points in $^{1}\mathcal{T}^{G}_{\lambda,\mu}$ whose fifth entry is strictly greater than zero. Hence we obtain the equality 
\begin{equation}\label{recursion2}|^{1}\mathcal{T}^{G}_{\lambda,\mu}|=|^{1}\mathcal{T}^{G}_{\lambda+\varpi_1-\varpi_2,\mu-\varpi_1}|+|^{2}\mathcal{T}^{G}_{\lambda,\mu}|,\end{equation}
where $^{2}\mathcal{T}^{G}_{\lambda,\mu}$ consists of tuples $(a,b,c,y) \in \mathbb{Z}_{+}^4$  (we substitute $y=m_2-b+d$) satisfying the following inequalities 
\begin{gather*}
	\label{frt1}	a+2b+c+y  \le n_1+m_2, \qquad c \le 1, \qquad a+2y\leq m_1+2m_2, \\
	\label{frt2}	c+2b+2y\leq m_1+2m_2,\qquad y+b\geq m_2.
	\end{gather*}	
In the next step we have to work a little bit more. We define once more an injective map 
$$^{2}\mathcal{T}^{G}_{\lambda-\varpi_1,\mu-\varpi_1}\rightarrow \ ^{2}\mathcal{T}^{G}_{\lambda,\mu},\ (a,b,c,y)\mapsto (a+1,b,c,y),$$
whose image consists exactly of those points in $^{2}\mathcal{T}^{G}_{\lambda,\mu}$ whose first entry is strictly greater than zero and $c+2b+2y<m_1+2m_2$. Hence we obtain
\begin{equation}\label{recursion3}|^{2}\mathcal{T}^{G}_{\lambda,\mu}|=|^{2}\mathcal{T}^{G}_{\lambda-\varpi_1,\mu-\varpi_1}|+|^{3}\mathcal{T}^{G}_{\lambda,\mu}|,\end{equation}
where $^{3}\mathcal{T}^{G}_{\lambda,\mu}$ consists of all elements in $^{2}\mathcal{T}^{G}_{\lambda,\mu}$  where $a=0$ or $c+2b+2y=m_1+2m_2$. 
A straightforward calculation shows that the cardinality of $^{3}\mathcal{T}^{G}_{\lambda,\mu}$ is given by the cardinality of the disjoint union 
\begin{align*}
	Y_{\lambda,\mu}^1&=\left\{(b,c,y)\in\mathbb{Z}_+^3: \substack{\displaystyle 2b+c+y\leq n_1+m_2,\ c\leq 1\\ \\ \displaystyle c+2b+2y\leq m_1+2m_2,\ y+b\geq m_2}\right\},
	\\[0,2cm]
	Y_{\lambda,\mu}^2& =
	\left\{(a,y)\in\mathbb{Z}_+^2:\substack{\displaystyle a-y< n_1-m_2-m_1\\ \\ \displaystyle  a+2y< m_1+2m_2}\right\}.
\end{align*}
Note that the polytope $Y_{\lambda,\mu}^1$ corresponds to all elements in $^{2}\mathcal{T}^{G}_{\lambda,\mu}$ where $a=0$ and $Y_{\lambda,\mu}^2$ to all elements where $a>0$ and $c+2b+2y=m_1+2m_2$.
We summarize the combinatorics of $\mathcal{T}^{G}_{\lambda,\mu}$ in the following lemma.
\begin{lem}\label{334g} Let $\lambda,\mu\in P^+$. Then 
$$|\mathcal{T}^{G}_{\lambda,\mu}|=|\mathcal{T}^{G}_{\lambda-\varpi_1,\mu-\varpi_1}|+|^{1}\mathcal{T}^{G}_{\lambda+\varpi_1-\varpi_2,\mu-\varpi_1}|+|^{2}\mathcal{T}^{G}_{\lambda-\varpi_1,\mu-\varpi_1}|+|^{3}\mathcal{T}^{G}_{\lambda,\mu}|.$$
\proof This follows from \eqref{recursion1}, \eqref{recursion2} and \eqref{recursion3}.
\endproof
\end{lem}
\subsection{} We want to prove a similar result as Lemma~\ref{334g} for the polytope $\mathcal{S}^{G}_{\lambda,\mu}$. First consider the injective map
$$\mathcal{S}^{G}_{\lambda-\varpi_1,\mu-\varpi_1}\rightarrow \mathcal{S}^{G}_{\lambda,\mu},\ (a,b,c,d,e,f)\mapsto (a+1,b,c,d,e+1,f),$$
whose image consists exactly of those points in $\mathcal{S}^{G}_{\lambda,\mu}$ satisfying $a,e>0$. Hence we obtain 
\begin{equation}\label{recursion12}|\mathcal{S}^{G}_{\lambda,\mu}|=|\mathcal{S}^{G}_{\lambda-\varpi_1,\mu-\varpi_1}|+|^{1}\mathcal{S}^{G}_{\lambda,\mu}|,\end{equation}
where $^{1}\mathcal{S}^{G}_{\lambda,\mu}= R^1_{\lambda,\mu}\cup\ R^2_{\lambda,\mu}$ with 
\begin{gather*}
	R^1_{\lambda,\mu}=\left\{(b,c,d,e)\in\mathbb{Z}_+^4: \substack{\displaystyle b+e\leq m_2,\ c+d\leq m_1+m_2,\ b+c\leq m_1+m_2\\ \\ \displaystyle b+c+d+e\leq n_1,\ 2c+3d-b\leq m_1+n_1,\\ \\ \displaystyle b+2c+d\leq m_1+n_1}\right\},
	\\[0,2cm]
	R^2_{\lambda,\mu}=\left\{(a,b,c,d)\in\mathbb{Z}_+^4: \substack{\displaystyle a<m_1,\ b\leq m_2,\ a+c+d< m_1+m_2\\ \\ \displaystyle a+b+c<m_1+m_2,\ a+b+c+d<n_1,\\ \\ \displaystyle 2a+2c+3d-b<m_1+n_1-1,\ 2a+b+2c+d<m_1+n_1-1}\right\}.
\end{gather*}

Note that the first polytope correponds to the elements in $^{1}\mathcal{S}^{G}_{\lambda,\mu}$ where $a=0$ and the second polytope to the elements where $a>0$ and $e=0$. We denote by $R^2_{\lambda,\mu}(a>0)$ the subset of elements with $a>0$. Again we define an injective map 
$$R^1_{\lambda+\varpi_1-\varpi_2,\mu-\varpi_1}\cup R^2_{\lambda+\varpi_1-\varpi_2,\mu-\varpi_1}(a>0)\rightarrow \ ^{1}\mathcal{S}^{G}_{\lambda,\mu}$$
where an element $(b,c,d,e)$ from $R^1_{\lambda+\varpi_1-\varpi_2,\mu-\varpi_1}$ is mapped to $(b,c,d,e+1)$ and an element $(a,b,c,d)$ from $R^2_{\lambda+\varpi_1-\varpi_2,\mu-\varpi_1}(a>0)$ is mapped to $(a-1,b+1,c,d+1)$. Hence the image consists exactly of those points in $^{1}\mathcal{S}^{G}_{\lambda,\mu}$ where $e>0$ in $R^1_{\lambda,\mu}$ and $b,d>0$ in $R^2_{\lambda,\mu}$. We get 
$$|^{1}\mathcal{S}^{G}_{\lambda,\mu}|=|^{1}\mathcal{S}_{\lambda+\varpi_1-\varpi_2,\mu-\varpi_1}|-|R^2_{\lambda+\varpi_1-\varpi_2,\mu-\varpi_1}(a=0)|+|R^1_{\lambda,\mu}(e=0)|+|R^2_{\lambda,\mu}(bd=0)|.$$
Obviously $R^1_{\lambda,\mu}(e=0)$ is described by the inequalities 
$$b\leq m_2,\ c+d\leq m_1+m_2,\ b+c\leq m_1+m_2,$$  
$$b+c+d\leq n_1,\ 2c+3d-b\leq m_1+n_1,\ b+2c+d\leq m_1+n_1.$$
Since
$$R^2_{\lambda+\varpi_1-\varpi_2,\mu-\varpi_1}(a=0)\rightarrow R^1_{\lambda,\mu}(e=0),\ (b,c,d)\mapsto (b+1,c,d+1)$$
is an injective map, whose image consists of all points in $R^1_{\lambda,\mu}(e=0)$ satisfying $b,d>0$, we are left to describe the disjoint union
$$^{2}\mathcal{S}^G_{\lambda,\mu}:=R^2_{\lambda,\mu}(bd=0)\cup R^1_{\lambda,\mu}(e=0,bd=0).$$
Writing down the inequalities it is straightforward to see that $^{2}\mathcal{S}^G_{\lambda,\mu}$ can be identified with the disjoint union of the following two polytopes
\begin{gather*}
	Q_{\lambda,\mu}^1=\left\{(a,c,d)\in\mathbb{Z}_+^3: \substack{\displaystyle a\leq m_1,\ a+c+d\leq \min\{m_1+m_2,n_1\}\\ \\ \displaystyle 2(a+c)+3d\leq m_1+n_1}\right\},
	\\[0,2cm]
	Q_{\lambda,\mu}^2=\left\{(a,b,c)\in\mathbb{Z}_+^3: \substack{\displaystyle a\leq m_1,\ a+b+c\leq \min\{m_1+m_2-1,n_1-1\}\\ \\ \displaystyle b\leq m_2-1,\ 2(a+c)+b\leq m_1+n_1-1}\right\}.
\end{gather*}
$Q_{\lambda, \mu}^1$ corresponds to $R^2_{\lambda, \mu}(b = 0) \cup R_{\lambda, \mu}^1(e = 0, b = 0)$ and 
$Q_{\lambda, \mu}^2$ to $R^2_{\lambda, \mu}(d = 0, b > 0)\cup R_{\lambda, \mu}^1(e =0, d = 0, b > 0)$.
We have
\begin{equation}\label{recursion22}|^{1}\mathcal{S}^{G}_{\lambda,\mu}|=|^{1}\mathcal{S}^{G}_{\lambda+\varpi_1-\varpi_2,\mu-\varpi_1}|+|^{2}\mathcal{S}^{G}_{\lambda,\mu}|.\end{equation}
In the next step we define again an injective map 
$$^{2}\mathcal{S}^{G}_{\lambda-\varpi_1,\mu-\varpi_1}\rightarrow \ ^{2}\mathcal{S}^{G}_{\lambda,\mu},$$
where an element $(a,c,d)$ from $Q_{\lambda-\varpi_1,\mu-\varpi_1}^1$ is mapped to $(a+1,c,d)$ and an element $(a,b,c)$ from $Q_{\lambda-\varpi_1,\mu-\varpi_1}^2$ is mapped to $(a+1,b,c)$. Obviously 
\begin{equation}\label{recursion32}|^{2}\mathcal{S}^{G}_{\lambda,\mu}|=|^{2}\mathcal{S}^{G}_{\lambda-\varpi_1,\mu-\varpi_1}|+|^{3}\mathcal{S}^{G}_{\lambda,\mu}|,\end{equation}
where $^{3}\mathcal{S}^{G}_{\lambda,\mu}$ consists of all elements in $^{2}\mathcal{S}^{G}_{\lambda,\mu}$  such that $a=0$. Hence it can be identified with the disjoint union of 
\begin{align*}
	\left\{(c,d)\in\mathbb{Z}_+^2: \substack{\displaystyle  c+d\leq \min\{m_1+m_2,n_1\}\\ \\ \displaystyle 2c+3d\leq m_1+n_1}\right\}, 
&& \left\{(b,c)\in\mathbb{Z}_+^2: \substack{\displaystyle b+c<\min\{m_1+m_2,n_1\}\\ \\ \displaystyle b< m_2,\ 2c+b< m_1+n_1}\right\}.
\end{align*}
We summarize the combinatorics of $\mathcal{S}^{G}_{\lambda,\mu}$ in the following lemma.
\begin{lem}\label{334g3} Let $\lambda,\mu\in P^+$. Then 
$$|\mathcal{S}^{G}_{\lambda,\mu}|=|\mathcal{S}^{G}_{\lambda-\varpi_1,\mu-\varpi_1}|+|^{1}\mathcal{S}^{G}_{\lambda+\varpi_1-\varpi_2,\mu-\varpi_1}|+|^{2}\mathcal{S}^{G}_{\lambda-\varpi_1,\mu-\varpi_1}|+|^{3}\mathcal{S}^{G}_{\lambda,\mu}|.$$
\begin{proof} This follows from \eqref{recursion12}, \eqref{recursion22} and \eqref{recursion32}.
\end{proof}
\end{lem}
\subsection{}The following lemma finishes the proof of Proposition~\ref{wichtp} in the case of $G_2$.
\begin{lem}\label{}
There exists a bijective map $\mathcal{S}^{G}_{\lambda,\mu} \rightarrow  \mathcal{T}^{G}_{\lambda,\mu}$.
\begin{proof}
In view of Lemma~\ref{334g}, Lemma~\ref{334g3} and an easy induction argument it will be enough to prove that $|^{3}\mathcal{S}^{G}_{\lambda,\mu}|$ and  $|^{3}\mathcal{T}^{G}_{\lambda,\mu}|$ have the same cardinality which we prove by induction on $m_2$.
\vspace{0,2cm}
\textit{Base step:} If $m_2=0$, then the cardinality of  $^{3}\mathcal{S}^{G}_{\lambda,\mu}$ is given by
\begin{equation}\label{3291}\sum_{i=0}^{M}\min\left\{M+1-i,\left\lfloor \frac{m_1+n_1-2i}{3}\right\rfloor+1\right\}\end{equation}
where $M=\min\{m_1,n_1\}$. 
If $M=0$ the statement is clear. $M=1$, then it is easy to see that $|^{3}\mathcal{T}^{G}_{\lambda,\mu}|=2$ if $m_1=n_1=1$ and $|^{3}\mathcal{T}^{G}_{\lambda,\mu}|=3$ otherwise. Since the same is true for \eqref{3291}, we obtain the claim for $M=1$. We assume from now on that $M\geq 2$ and define an injective map 
$$^{3}\mathcal{T}^{G}_{\lambda-2\varpi_1,\mu-\varpi_1}\rightarrow ^{3}\mathcal{T}^{G}_{\lambda,\mu},$$
where an element $(b,c,y)$ from $Y_{\lambda-2\varpi_1,\mu-\varpi_1}^1$ is mapped to $(b,c,y+1)$ and an element $(a,y)$ from $Y_{\lambda-2\varpi_1,\mu-\varpi_1}^2$ is mapped to $(a,y+1)$. We have 
$$|^{3}\mathcal{T}^{G}_{\lambda,\mu}|=|^{3}\mathcal{T}^{G}_{\lambda-2\varpi_1,\mu-\varpi_1}|+\min\{n_1,2m_1\}+1,$$
where the last summand counts the elements in $^{3}\mathcal{T}^{G}_{\lambda,\mu}$ with $y=0$. By induction we have 
\begin{align*}|^{3}\mathcal{T}^{G}_{\lambda,\mu}|&=\sum_{i=0}^{\min\{m_1-1,n_1\}-1}\min\left\{\min\{m_1-1,n_1\}-i,\left\lfloor \frac{m_1+n_1-2i}{3}\right\rfloor\right\}+\min\{n_1,2m_1\}+1&\\&=\sum_{i=0}^{\min\{m_1-1,n_1\}}\min\left\{\min\{m_1-1,n_1\}-i,\left\lfloor \frac{m_1+n_1-2i}{3}\right\rfloor\right\}+\min\{n_1,2m_1\}+1&\\&=\begin{cases}\displaystyle\sum_{i=0}^{M}\min\left\{M-i+1,\left\lfloor \frac{m_1+n_1-2i}{3}\right\rfloor+1\right\},& \text{if $n_1<m_1$},\\
\displaystyle\sum_{i=0}^{m_1}\min\left\{m_1-i,\left\lfloor \frac{m_1+n_1-2i}{3}\right\rfloor+1\right\}+n_1-m_1+1,& \text{if $m_1\leq  n_1<2m_1$},\\
\displaystyle\sum_{i=0}^{m_1}\min\left\{m_1-i,\left\lfloor \frac{m_1+n_1-2i}{3}\right\rfloor+1\right\}+m_1+1,& \text{if $2m_1\leq n_1$}
\end{cases}&\\&=\sum_{i=0}^{M}\min\left\{M+1-i,\left\lfloor \frac{m_1+n_1-2i}{3}\right\rfloor+1\right\},\end{align*}
where the last equality follows in the last two cases from the fact that 
$m_1-i < \left\lfloor \frac{m_1 + n_1 - 2i}{3}\right\rfloor+1$ if and only if $2m_1-n_1\leq i\leq m_1$. This finishes the base step.

\vspace{0,2cm}
\textit{Inductive step:}
If $m_2>0$ we have
$|^{3}\mathcal{S}^{G}_{\lambda,\mu}|=|^{3}\mathcal{S}^{G}_{\lambda-\varpi_2,\mu}|+k$
where $k$ is the cardinality of the union of the following sets
\begin{align*}
	Z_{\lambda, \mu}^1&= \left\{(c,d)\in\mathbb{Z}^2_+: \substack{\displaystyle 2c+3d\leq m_1+n_1,\\ \\ \displaystyle \min\{m_1+m_2-1,n_1\}<c+d\leq \min\{m_1+m_2,n_1\}}\right\},	\\[0,2cm] 
	Z_{\lambda, \mu}^2&=\left\{(b,c)\in\mathbb{Z}_+^2: \substack{\displaystyle 2c+b< m_1+n_1,\ \ b<m_2-1,\\ \\ \displaystyle \min\{m_1+m_2-1,n_1\}\leq b+c<\min\{m_1+m_2,n_1\}}\right\},\\[0,2cm] 
	Z_{\lambda, \mu}^3&=\left\{c\in\mathbb{Z}_+:  c\leq \min\{m_1,n_1-m_2\}\right\}.
\end{align*}
The cardinalities of these sets are given by 
\begin{align*}
	|Z_{\lambda, \mu}^1| & = 
	\begin{cases}
		0,&\text{if $n_1<m_1+2m_2$},\\
		n_1-m_1-2m_2+1,& \text{if $m_1+2m_2\leq n_1<2m_1+3m_2$},\\
		m_1+m_2+1,& \text{if $ 2m_1+3m_2\leq n_1$}.
	\end{cases}\\ 
	|Z_{\lambda, \mu}^2| &= 
	\begin{cases}
		0,&\text{if $n_1<m_1+m_2$},\\
		n_1-m_1-m_2,& \text{if $m_1+m_2\leq n_1<m_1+2m_2$},\\
		m_2-1,& \text{if $ m_1+2m_2\leq n_1$} \\ 
	\end{cases}\\
	|Z_{\lambda, \mu}^3| &= 
	\begin{cases}
		0,&\text{if $n_1< m_2$},\\
		(\min\{m_1,n_1-m_2\}+1), & \text{if $n_1\geq m_2$}.
	\end{cases}
\end{align*}
So adding these numbers up, we get that 
\begin{align*}
	k = |Z_{\lambda, \mu}^1| + |Z_{\lambda, \mu}^2| + |Z_{\lambda, \mu}^3| = 
	\begin{cases}
		0,&\text{if $n_1<m_2$},\\
		n_1-m_2+1,&\text{if $m_2\leq n_1<2m_1+3m_2$},\\
		2(m_1+m_2)+1,& \text{if $ 2m_1+3m_2\leq n_1$}.
	\end{cases}
\end{align*}
On the other hand we have an injective map 
$$^{3}\mathcal{T}^{G}_{\lambda-\varpi_2,\mu}\rightarrow ^{3}\mathcal{T}^{G}_{\lambda,\mu},$$
where an element $(a,y)$ from $Y_{\lambda-\varpi_2,\mu}^2$ is mapped to $(a,y+1)$ and an element $(b,c,y)$ from $Y_{\lambda-\varpi_2,\mu}^1$ is mapped to $(b,c,y+1)$. So we also have 
$|^{3}\mathcal{T}^{G}_{\lambda,\mu}|=|^{3}\mathcal{T}^{G}_{\lambda-\varpi_2,\mu}|+m$
where $m$ is the cardinality of the union of the following sets:
\begin{gather*}
	\left\{(b,c)\in\mathbb{Z}^2_+: 2b+c\leq \min\{m_2+n_1,m_1+2m_2\},\ \ c\leq 1,\  \ b\geq m_2\right\},\\ 
	\left\{a\in\mathbb{Z}_+: a< \min\{n_1-m_2-m_1,m_1+2m_2\}\right\}.
\end{gather*}
The cardinality of the first set is $0$ if $n_1<m_2$ and otherwise $\min\{n_1-m_2,m_1\}+1$. The cardinality of the second set is given by $0$ if $n_1< m_1+m_2$ and otherwise by $\min\{n_1-m_2-m_1,m_1+2m_2\}.$ This shows $k=m$ and hence the claim. 
\end{proof}
\end{lem}
\bibliographystyle{plain}
\bibliography{bibfile}

\end{document}